\pdfoutput=1

\documentclass[12pt,longtitle,times]{amsart}

\usepackage{psfrag}
\usepackage{url}
\usepackage{caption}
\usepackage{subcaption}
\usepackage{parskip}
\usepackage{microtype}
\usepackage{amsfonts}
\usepackage{fleqn}
\usepackage{geometry}
\usepackage{graphicx}
\usepackage{graphicx,epstopdf}
\usepackage{amsmath}
\usepackage{amssymb,amsthm}
\usepackage{mathrsfs}
\usepackage{psfrag}
\usepackage[resetlabels]{multibib}
\usepackage{textcase}
\usepackage{txfonts}
\usepackage{xspace}
\usepackage{xcolor}
\usepackage{lipsum}
\usepackage[bottom,norule,hang]{footmisc}
\usepackage{setspace}

\newcommand{\Mod}[1]{\ (\text{mod}\ #1)}
\makeatletter
\let\@wraptoccontribs\wraptoccontribs
\makeatother

\begin{document}


\title[A prime fractal]{A prime fractal and global quasi-self-similar structure in the distribution of prime-indexed primes}

\author[R.G.~Batchko]{Robert G. Batchko\\
\\
\textit{H\lowercase{olochip} C\lowercase{orporation}, 4940 W. 147\lowercase{th} S\lowercase{treet}, H\lowercase{awthorne}, C\lowercase{alifornia} 90250 USA}\\
\lowercase{\texttt{rgb@holochip.com}}\\
}

\thanks{The author is indebted to Charles Kinzer for valuable discussions on statistics, and Jei-Yin Yiu and Samuel Robinson for helpful discussions on programming math software.}

\begin{abstract}
Let $p_n$ be the $n$th prime and $p_{p_n}$ be the $n$th prime-indexed prime (PIP). The process of taking prime-indexed subsequences of primes can be iterated, and the number of such iterations is the prime-index order. We report empirical evidence that the set composed of finite-differenced PIP sequences of prime-index order $k\geq{}1$ forms a quasi-self-similar fractal structure with scaling by prime-index order. Strong positive linear correlation ($r\geq{}0.926$) is observed for all pairwise combinations of these finite-differenced PIP sequences over the range of our sample, the first 1.3 billion primes. The structure exhibits translation invariance for shifts in the index set of the PIP sequences. Other free parameters of the structure include prime-index order and the order and spacing of the finite difference operator.  The structure is graphed using 8-bit color fractal plots, scaled across prime-index orders $k=1..6$ and spans the first 1.3 billion primes. 
\end{abstract}
\maketitle

\section{Overview}\label{Section:Introduction}

The search for structure amidst randomness in the distribution of the primes has been a quest of mathematicians since antiquity. To date, three forms of structure are known to exist. The first, \textit{local structure} (c. 200 B.C.), generally refers to information that is localized on some parameter of a given relation, rule or pattern, and includes residue classes\cite{Granville-Different-Approaches-2009,Tao_The-Dichotomy-ICM2005} and arithmetic progressions\cite{Green-ARbitr-long-ARITH-2008,Tao-ARbitr-long-poly-2008}. For example, in an Eratosthenes sieve, the step of eliminating even numbers is equivalent to eliminating ``local information at the 2 place" -- T. Tao.\footnote{\label{Tao-structure-randomness2009}T. Tao, Structure and randomness in the prime numbers, (2009), Blog post on https://terrytao.wordpress.com/2009/07/20/structure-and-randomness-in-the-prime-numbers-2/.}

The second form is \textit{asymptotic (or large-scale) structure}, the key result of the proof of the Prime Number Theorem. Asymptotic structure refers to the fact that, as the primes grow in size, the gaps between them become increasingly large and their density becomes increasingly structured.\cite{Goldstein-1973}\footnote{Asymptotic structure can also be considered ``local information at the infinity place'' -- T. Tao, see footnote \ref{Tao-structure-randomness2009}.}

The third and most recent form we will call \textit{statistical structure}, wherein authors have reported empirical evidence of fractal\cite{Cattani-Fractal-2010}, chaotic\cite{Bershadskii-Hidden-2011,Wolf1996}, and other non-random\cite{Ares-Hidden-2006,Holdom-Scale-Invariant-2009,Selvam-Universal-2014,Szpiro-TheGaps-2004,Szpiro-Peaks-and_Gaps-2007}\footnote{\label{Liang2006}W. Liang and H. Yan, arXiv:math/0603450v1 [math.NT] (2006).}$^,$\footnote{\label{Timberlake2008}T. K. Timberlake and J. M. Tucker, arXiv:0708.2567v2 [quant-ph] (2008).} behavior in the distribution of the primes. However, due to the absence of an actual source behind any of those data (e.g., some non-obvious structure or pattern in the primes themselves), that work was generally limited to examining only the probability distributions of functions acting on the primes.\cite{Bershadskii-Hidden-2011,Dahmen-Similarity-In-The-Statistics-2001,Szpiro-TheGaps-2004,Szpiro-Peaks-and_Gaps-2007} As a result of this paradox -- substantial indirect evidence of structure but no sighting of the structured object itself -- statistical structure has been generally regarded more as a mathematical curiosity than an underlying property of the primes.\footnote{\label{Cloitre2011}B. Cloitre, On the fractal behavior of primes, (2011), Available from http://bcmathematics.monsite-orange.fr/FractalOrderOfPrimes.pdf.}$^,$\footnote{The notion that a hidden fractal structure may be embedded in the distribution of the primes has received increased attention since Folsom, Kent and Ono (2012) showed that the values of the partition function $p(n)$ are $\ell$-adically fractal for primes $\ell \geq 5$.\cite{Bruinier-Algebraic-2013,Folsom-L-Adic-2012}}

In addition to the above-mentioned three existent forms of structure --  \textit{local}, \textit{asymptotic}, and \textit{statistical} -- a fourth, mythical, form was offered by Tao in 2006. He describes \textit{exotic structure} as a hypothetical property of the primes, in which, ``they [the primes] could obey some exotic structure not predicted by [the Cram\'{e}r random] model, e.g., they could be unexpectedly dense on some structured set,'' while noting, ``[Because] We don't know if they also have some additional exotic structure \dots, we have been unable to settle many questions about primes.''\cite{Pintz-Cramer-vs-Cramer-2007,Tao-Long-2006,Tao-Recent-Progress-2009}

In this manuscript, we report the first empirical observation of a possible (quasi-) exotic structure and source of statistical structure -- a global pseudorandom fractal object encoded in distribution of the primes -- located at the crossroads of the prime-indexed primes and finite differences.

The remainder of this manuscript is organized as follows. Section 2 gives a review of previous work and introduces the family of finite-differenced prime-indexed-prime sequences $\daleth$ and its quasi-self-similar behavior. In Section 3, using a sample of the first $1.3\times{}10^9$ primes, statistical analysis is performed on $\daleth$, and strong positive correlation is observed among all pairwise combinations of its sequences. In Section 4, the fractal nature of $\daleth$ is presented graphically, both in its raw form and after transformation to stationarize its variance. In Section 5, the zeros, free parameters and global nature of $\daleth$ are examined before concluding.

\section{Background and Introduction}
\subsection{Prime-Indexed Primes (PIPs).} Following Dressler and Parker\cite{Dressler1975}, Broughan and Barnett\cite{Broughan2009}, and others\cite{Bayless2013,Jordan-On-Sums-1965}, the prime-indexed primes (PIPs) $q_i$ are defined as follows, 
\begin{center}
``If $p_i$ is the $i$th prime, then we define $q_i$ to be $p_{p_i}$.''
 \end{center}
 
  Thus, for $i = 1,2,3,\dots$, we have $q_i = \left(p_{p_1},p_{p_2},p_{p_3},\dots\right) = \left(p_2,p_3,p_5,\dots\right) = \left(3,5,11,\dots\right)$.\footnote{To date, the terminology for the set $q_i$ has been somewhat non-standardized; besides ``prime-indexed primes'' (``PIPs'') used here and previously by Broughan and Barnett\cite{Broughan2009} and Bayless et al.,\cite{Bayless2013} it has been called ``superprimes'', ``higher-order primes'' and ``primeth primes'', just to name a few (see: The On-Line Encyclopedia of Integer Sequences, 2012, \url{http://oeis.org/wiki/Higher-order_prime_numbers}[Online; accessed October 01, 2012] and Wikipedia, Super-prime, 2013, \url{http://en.wikipedia.org/wiki/Super-prime, [Online; accessed  13-December-2013]}.}
  
Several authors have noted that, ``the process of taking a prime indexed subsequence [of the primes] can be iterated.''\cite{Broughan2009} For example, taking the prime-indexed subsequence of $q_i$ gives $q_{p_i} = \left(p_{p_{p_1}},p_{p_{p_2}},p_{p_{p_3}},\dots\right) = \left(p_{p_2},p_{p_3},p_{p_5},\dots\right) = \left(p_3,p_5,p_{11},\dots\right) = \left(5,11,31,\dots\right)$. Let us call the number of such iterations the \textit{prime-index order}.\footnote{Fernandez defined an ``order of primeness'' for the PIPs, which is not to be confused with the definition of prime-index order used here (Fernandez, N., \textit{An order of primeness, f(p)}, unpublished (1999), available at http://borve.org/primeness/FOP.html).}

Here, we are interested in addressing the value of $q_i$ for any prime-index order. Thus, we redefine $q_i$ with the augmentation of two additional parameters as follows:

\begin{singlespace}
\begin{equation}
\centering
\hphantom{------------}q_{si}^k = q_s^{k}(i)=
\begin{cases}\label{eq:Cases}
    i=1,2,\dots        & \text{if } k = 0,\\
    p_i,             & \text{if } k = 1,\\
    p_{p_i},         & \text{if } k = 2,\\
    \vdots          & \vdots\\
    p_{p_{\ddots_{p_i}}},   & \text{if } k > 2                
\end{cases} 
\end{equation}
\end{singlespace}

where $i = \{1,2,\dots$\}, $k$ is the prime-index order ($k \in \mathbb{N}_0$), and $s$, which will be discussed in Section \ref{Section:Discussion}, is a shift parameter of the index set of $q_s^k(i)$. Accordingly, future references to ``PIPs'' will follow our generalized definition of $q$.

Table \ref{tab:q_k^i} shows an array of the values of $q_0^k(i)$ for a domain $i = 1..20$ and $k = 0..8$.\footnote{Note that $q_{0}^0(i)$ and $q_{0}^1(i)$ are the sequences of positive integers and primes, respectively.}

\begin{table}[h]
\centering
\fontsize{10pt}{0}
\caption{Values of $q_0^k(i)$ for $i = 1..20$ and $k = 0..8$}
$\begin{array}{ccccccccc} k=0 & k=1 & k=2 & k=3 & k=4 & k=5 & k=6 & k=7 & k=8\\ 
\hline
1 & 2 & 3 & 5 & 11 & 31 & 127 & 709 & 5381\\ 2 & 3 & 5 & 11 & 31 & 127 & 709 & 5381 & 52711\\ 3 & 5 & 11 & 31 & 127 & 709 & 5381 & 52711 & 648391\\ 4 & 7 & 17 & 59 & 277 & 1787 & 15299 & 167449 & 2269733\\ 5 & 11 & 31 & 127 & 709 & 5381 & 52711 & 648391 & 9737333\\ 6 & 13 & 41 & 179 & 1063 & 8527 & 87803 & 1128889 & 17624813\\ 7 & 17 & 59 & 277 & 1787 & 15299 & 167449 & 2269733 & 37139213\\ 8 & 19 & 67 & 331 & 2221 & 19577 & 219613 & 3042161 & 50728129\\ 9 & 23 & 83 & 431 & 3001 & 27457 & 318211 & 4535189 & 77557187\\ 10 & 29 & 109 & 599 & 4397 & 42043 & 506683 & 7474967 & 131807699\\ 11 & 31 & 127 & 709 & 5381 & 52711 & 648391 & 9737333 & 174440041\\ 12 & 37 & 157 & 919 & 7193 & 72727 & 919913 & 14161729 & 259336153\\ 13 & 41 & 179 & 1063 & 8527 & 87803 & 1128889 & 17624813 & 326851121\\ 14 & 43 & 191 & 1153 & 9319 & 96797 & 1254739 & 19734581 & 368345293\\ 15 & 47 & 211 & 1297 & 10631 & 112129 & 1471343 & 23391799 & 440817757\\ 16 & 53 & 241 & 1523 & 12763 & 137077 & 1828669 & 29499439 & 563167303\\ 17 & 59 & 277 & 1787 & 15299 & 167449 & 2269733 & 37139213 & 718064159\\ 18 & 61 & 283 & 1847 & 15823 & 173867 & 2364361 & 38790341 & 751783477\\ 19 & 67 & 331 & 2221 & 19577 & 219613 & 3042161 & 50728129 & 997525853\\ 20 & 71 & 353 & 2381 & 21179 & 239489 & 3338989 & 56011909 & 1107276647 \end{array}$
\label{tab:q_k^i}
\end{table}

\subsubsection[short]{Asymptotic Behavior of PIPs.} Building on the PNT, progress has been made in deriving the asymptotic form for PIPs where prime-index order $k=2$.\footnote{Above, we introduced $s$ as a shift parameter for the index set of PIPs, which will be discussed in Section \ref{Section:Discussion}. For now, $s=0$, indicating no shift.} Broughan and Barnett \cite{Broughan2009} derived an upper bound of 

\begin{equation}\label{eq:Eqn-Broughan}
\hphantom{--------}q_{0}^2(n)=q_{0,n}^2 = n\text{log}^2n + 3n\text{log}n\text{log}\text{log}n+O(n\text{log}n) \sim n\text{log}^2n.
\end{equation}

By iteration of \eqref{eq:Eqn-Broughan}, the asymptotic behavior of $q_{0,n}^k$ for any prime-index order $k$ can be shown to be\footnote{The On-Line Encyclopedia of Integer Sequences, 2013, \url{http://oeis.org/wiki/Higher-order_prime_numbers}[Online; accessed October 01, 2012].}

\begin{equation}\label{eq:Eqn-Wiki-asymp-q}
\hphantom{----------------}q_{0,n}^k = \underbrace{p_{p_{\ddots_{p_{p_n}}}}}_{\text{$p$'s ($k$ times)}} \sim{} n(\text{log\space}n)^k.
\end{equation}
   
Refining this approximation, Bayless, Klyve and Silva\cite{Bayless2013} took Dusart's lower bound for $p_n$ (PIPs of $k=1$)

\begin{equation}\label{eq:Bayless-iteration-1}
\hphantom{---------------}q_{0,n}^1 = p_n > n(\text{log}n + \text{log}\text{log}n-1) 
\end{equation}

and iterated it to derive an improved lower bound for PIPs of $k=2$  
 
\begin{singlespace} 
\begin{eqnarray}\label{eq:Bayless-iteration-2}
\hphantom{---------}q_{0,n}^2 > n(\text{log}n + \text{log}\text{log}n-1)\left(\vphantom{p\ddots_p\ddots}\text{log}(n(\text{log}n + \text{log}\text{log}n-1) )\right. \\
\left.+\text{ }\vphantom{p\ddots_p\ddots}\text{log}\text{log}(n(\text{log}n + \text{log}\text{log}n-1))-1\right).\nonumber
\end{eqnarray}
\end{singlespace}

\begin{figure}[h]
\centering
\includegraphics[width=0.7\linewidth]{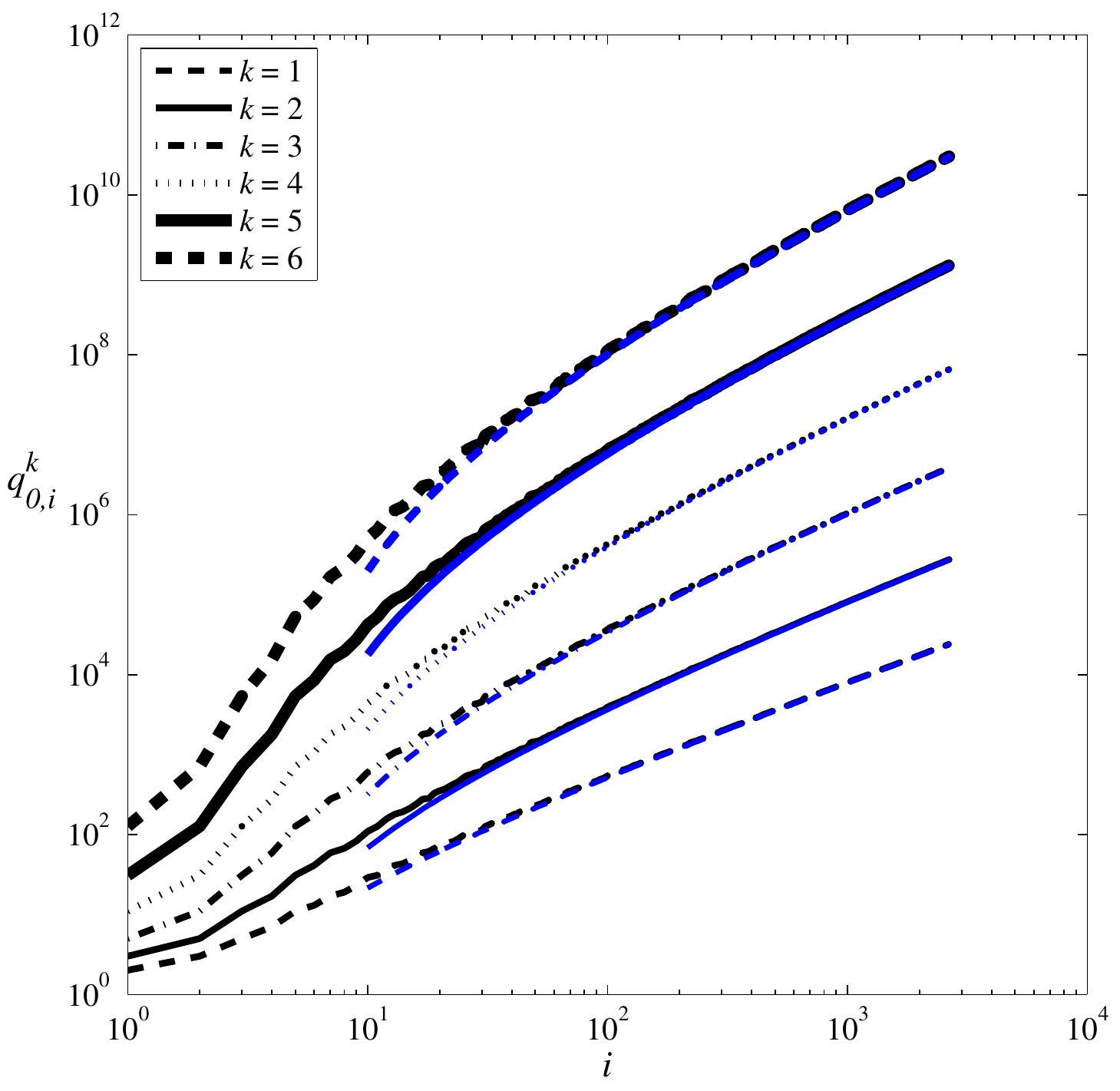}
\caption{Log-log plots of $q_{0,i}^k$ vs. $i$ for $i = 1..2500$ and $k = 1..6$}
\label{fig:Q_i_k_004}
\end{figure}

This iterative process can be continued to find the lower bound for $q_{0,n}^k$ for $k$ in general.\footnote{For each iteration, $n$ is replaced with $n(\text{log}n + \text{log}\text{log}n-1)$.} Figure \ref{fig:Q_i_k_004} contains six log-log plots showing the growth of $q_{0,i}^k$ vs. $i$ for $i = 1..2500$ and $k = 1..6$ (black lines). The asymptotic forms, based on iterations of \eqref{eq:Bayless-iteration-1} and \eqref{eq:Bayless-iteration-2}, for all six orders of $k$ are also shown (blue lines). 

\subsection[short]{Finite differences.} The forward finite difference operator $\Delta$, is the discrete analog to differentials of continuous functions.\footnote{In this discussion, we use forward differences to be consistent with other authors, however, backward or central differences could also be used.} In time series analysis, the forward first-order difference of a function $f(x)$ is defined as

\begin{equation*}
\Delta_h^1[f](x)=f(x+h)-f(x)
\end{equation*}

 where $h$ is the spacing parameter.\footnote{\label{Wiki-finite-difference}Wikipedia, Finite difference, 2014, \url{http://en.wikipedia.org/wiki/Finite_difference} [Online; accessed April 12, 2014].} First-order differences are commonly used to detrend series that are linear-nonstationary (e.g., non-stationary in the mean). Operation of $\Delta_1^1$ on sequences of consecutive primes 
 
 \begin{equation*}
 \Delta_1^1(p_i)=p_{i+1}-p_i,
 \end{equation*} 
 
 yields sequences of prime differences (or gaps), and is thus quintessential in the study of the distribution of primes.\cite{Erdos-The-Difference-1940,Erdos-Some-New-Questions-1948,Odlyzko-Jumping-1999,Odlyzyko-Iterated-Abs-Val-1993}
 
 Likewise, the forward second-order finite difference of a function $f(x)$
 
 \begin{equation*}
\Delta_h^2 [f](x)=f(x+2h)-2f(x+h)+f(x)
 \end{equation*}
 
 is used to detrend quadratic-nonstationary series. When operating on consecutive primes, the second-order finite difference operator
 
 \begin{equation*}
 \Delta_1^2(p_i)=p_{i+2}-2p_{i+1}+p_i
 \end{equation*}
 
 generates the sequence of prime increments (differences between consecutive prime gaps).
  
 In general, the $n$th-order forward finite difference of $f(x)$ is defined as 
 
 \begin{equation*}
 \Delta_h^n[f](x)=\sum_{m=0}^{n}(-1)^m\binom{n}{m}f(x+(n-m)h)
 \end{equation*}
 
where $n$ is the order of the finite difference operator and $\binom{n}{m}$ are the binomial coefficients.\footnote{See footnote \ref{Wiki-finite-difference}.} Szpiro\cite{Szpiro-TheGaps-2004,Szpiro-Peaks-and_Gaps-2007} examined the operation of $\Delta_1^n$ on the primes.
 
\subsubsection[short]{Finite differences of PIPs.} In this work, we are interested in the family of sequences obtained by taking the finite differences of PIPs.\footnote{Labos, E., Sequence A073131, ``a(n)=p[p[n+1]]-p[p[n]]] where p(j) is the j-th prime,'' in The On-Line Encyclopedia of Integer Sequences (2002), published electronically at http://oeis.org/A073131.} We thus define a \textit{general, $n$th-order prime-indexed-prime finite difference function} $\daleth$ as follows,

\begin{equation}\label{Defn daleth}
\hphantom{----------}\daleth_{hsi}^{nk}=\Delta_h^n(q_{si}^k)=\Delta_h^n[q_{s}^k](i)=\sum_{m=0}^{n}(-1)^m\binom{n}{m}q_{s}^k(i+(n-m)h),
\end{equation}

where $q_{si}^k$ is defined in \eqref{eq:Cases}. Table \ref{tab:Daleth_i_k} shows the range of

\begin{equation*}
\hphantom{----}\daleth_{1,0,i}^{2,k}=q_{0,i+2}^k-2q_{0,i+1}^k+q_{0,i}^k
\end{equation*}

for prime-index orders $k=0..8$, a domain of $i=1..50{}$, spacing of the finite difference operator $h=1$, second-order finite difference operator $n=2$, and zero shift ($s=0$) of the index set of $q_{0,i}^k$.

\begin{table}[b]
\centering
\caption{Values of $\daleth_{1,0,i}^{2,k}$ with $i = 1..50$ and $k = 0..8$}
\includegraphics[width=1\linewidth]{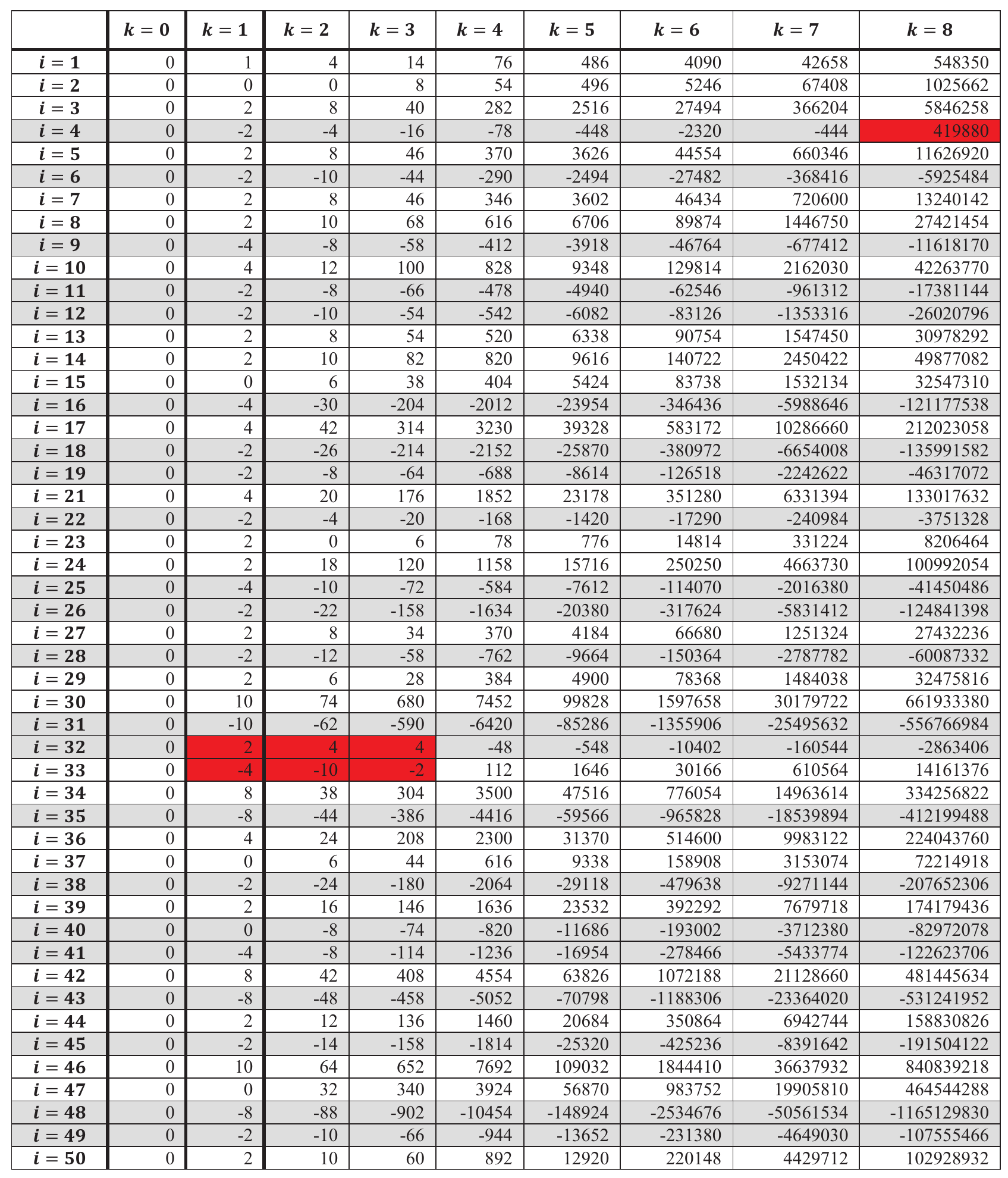}
\label{tab:Daleth_i_k}
\end{table}

In Table \ref{tab:Daleth_i_k}, it can be seen that in any given row $i$, all $k=1..8$ terms generally have the identical sign.\footnote{In the following discussion, we disregard the k=0 column of Table \ref{tab:Daleth_i_k}, in which all $i$ terms are zero due to the finite differencing of the integers.} This observation is highlighted by selectively shading the cells in Table \ref{tab:Daleth_i_k} as follows: the rows with mostly negative values are shaded in gray; rows with predominantly positive values are not shaded; and outliers (a term whose sign deviates from a majority of the others in the same row) are shaded in red.

\clearpage

The occurrence of only seven outliers in the 400 samples given in Table \ref{tab:Daleth_i_k} is unexpected. Given the pseudorandom nature of the distribution of the primes, one would anticipate approximately half (i.e., roughly 200) of the samples to be outliers. 

\subsubsection[short]{Time series plots of $\daleth_{1,0,i}^{2,k}$.} 
\begin{figure}[hp]
\centering
\includegraphics[width=1\linewidth]{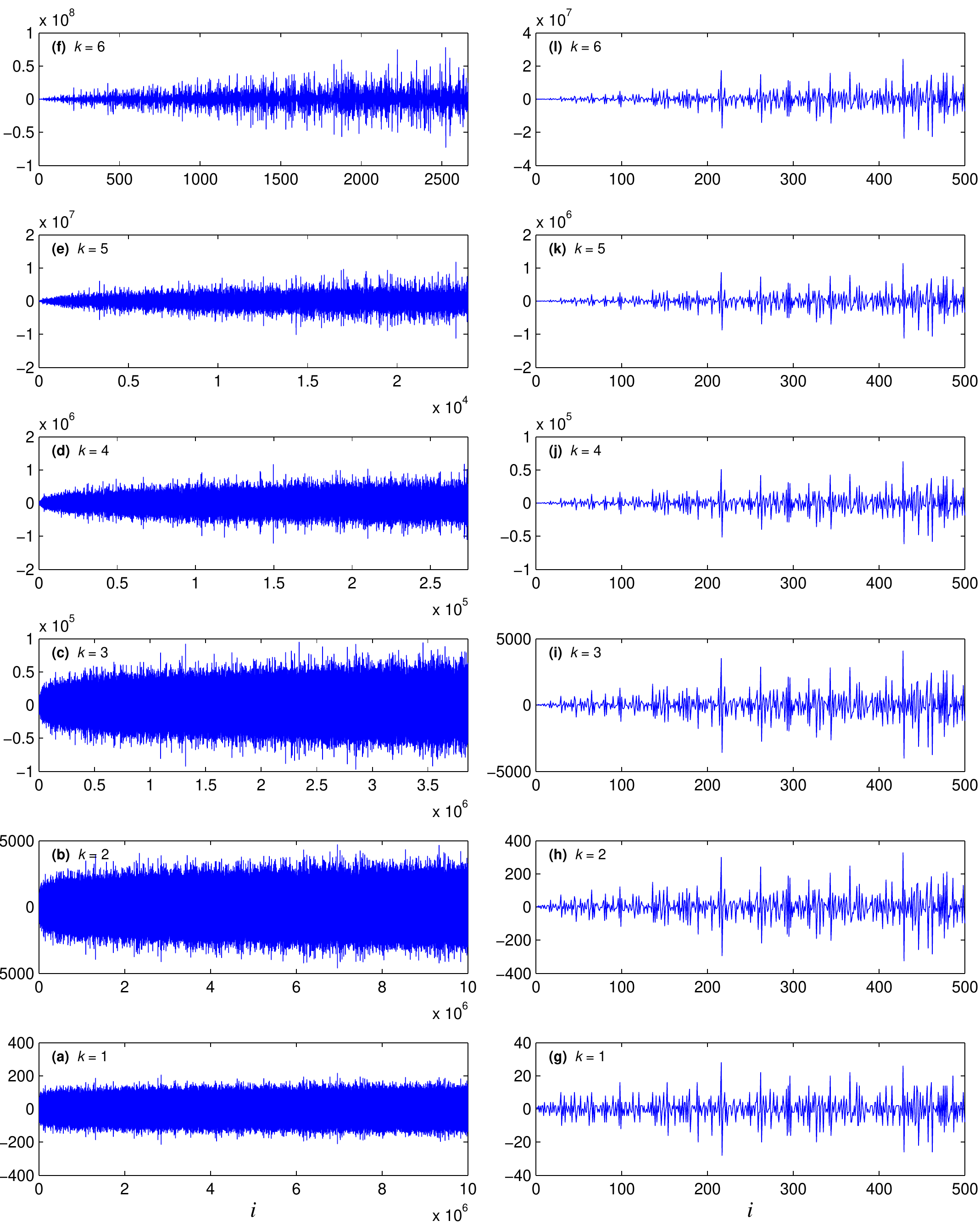}
\caption{$\daleth_{1,0,i}^{2,k}$ with $k=1..6$; (a-f) large samples; (g-l) sample size $T = 500$}
\label{fig:DELTA_Q_001}
\end{figure}

Expanding on the data given in Table \ref{tab:Daleth_i_k}, Figure \ref{fig:DELTA_Q_001} contains 12 plots of the evolution of $\daleth_{1,0,i}^{2,k}$ vs. $i$. Subplots (a)-(f) show $\daleth_{1,0,i}^{2,k}$ with $k=1..6$, respectively, and a domain of the full sample, $i=1..1.3\times{}10^9$.\footnote{Due to available processing power, we reduce the domain of $\daleth_{1,0,i}^{2,1}$ and $\daleth_{1,0,i}^{2,2}$ to $i=1..10^7.$} Due to the large number of data points and small page size, it is not possible to see fine detail in these plots; however, it is clear that the sequences rapidly oscillate around zero, slowly grow in variance over time, and have approximate mirror symmetry about the $\daleth_{1,0,i}^{2,k} = 0$ axis. To help visualize some finer detail of the sequences, subplots (g)-(l) show $\daleth_{1,0,i}^{2,k}$ with $k=1..6$, respectively, and a domain $i = 1..500$. It can be seen in subplots (g)-(l), that while $\daleth_{1,0,i}^{2,k}$ varies by six orders of magnitude as a function of $k$, all six sequences exhibit rapid oscillation about the $\daleth=0$ axis and follow an apparently similar structure in their growth and fine detail.

\clearpage

The seeming anomaly noted in Table \ref{tab:Daleth_i_k} and Figure \ref{fig:DELTA_Q_001} (unexpectedly similar patterns among the $k$th sequences of $\daleth_{1,0,i}^{2,k}$) inspires us to investigate whether \textit{the relationships between the distributions of the sequences of $\daleth_{1,0,i}^{2,k}$, and more generally those of $\daleth_{hsi}^{nk}$, are not the result of randomness}. 

\section{Statistical Analysis of $\daleth$}
We now examine statistical properties of the example case $\daleth_{1,0,i}^{2,k}$, where the finite difference spacing parameter $h=1$, order of finite differencing $n=2$, and PIP index-set shift parameter $s=0$. The domain in the following study is generally maximized based on the size of our sample, which consists of $p_1..p_{1.3\times{}10^9}$, where $p$ is prime. However, in some cases we reduce the domain due to limitations of processing power or display resolution, and in other cases we crop the domain appropriately so that all displayed sequences have the identical number of elements.

\subsection[short]{Variance and mean.} In developing a statistical model, it is of interest to determine whether the mean and variance of the distribution of the time series under study are stationary over the sample. To assess stationarity of these moments of $\daleth_{1,0,i}^{2,k}$, let $T$ be the sample size, $w$ be the width of a window (or sub-sample), $y$ be the step size and define the rolling-sample means $\hat{\mu}_i(w)$ and variances $\hat{\sigma}_i^2(w)$ as follows,

\begin{eqnarray}
\hphantom{------------}\nonumber&&\hat{\mu}_i(w) = \frac{1}{w}\sum_{j=0}^{w-1}\daleth_{1,i-j}^{2,k}\bigskip,\\
\nonumber&&\hat{\sigma}_i^2(w) = \frac{1}{w-1}\sum_{j=0}^{w-1}(\daleth_{1,i-j}^{2,k}-\hat{\mu}_i(w))^2,
\end{eqnarray}

for windows $i=w,w+y,\dots,w+\left\lfloor{}\frac{T-w}{y}\right\rfloor{}y$, where $\lfloor{}x\rfloor{}$ is the floor function. These $\hat{\mu}_i(w)$ and $\hat{\sigma}_i^2(w)$ are estimated moments based on the most recent $w$ observations, taken at time $i$, with a window width $w$, and stepped by $y$ samples at a time.

\begin{figure}[hp]
\centering
\includegraphics[width=1\linewidth]{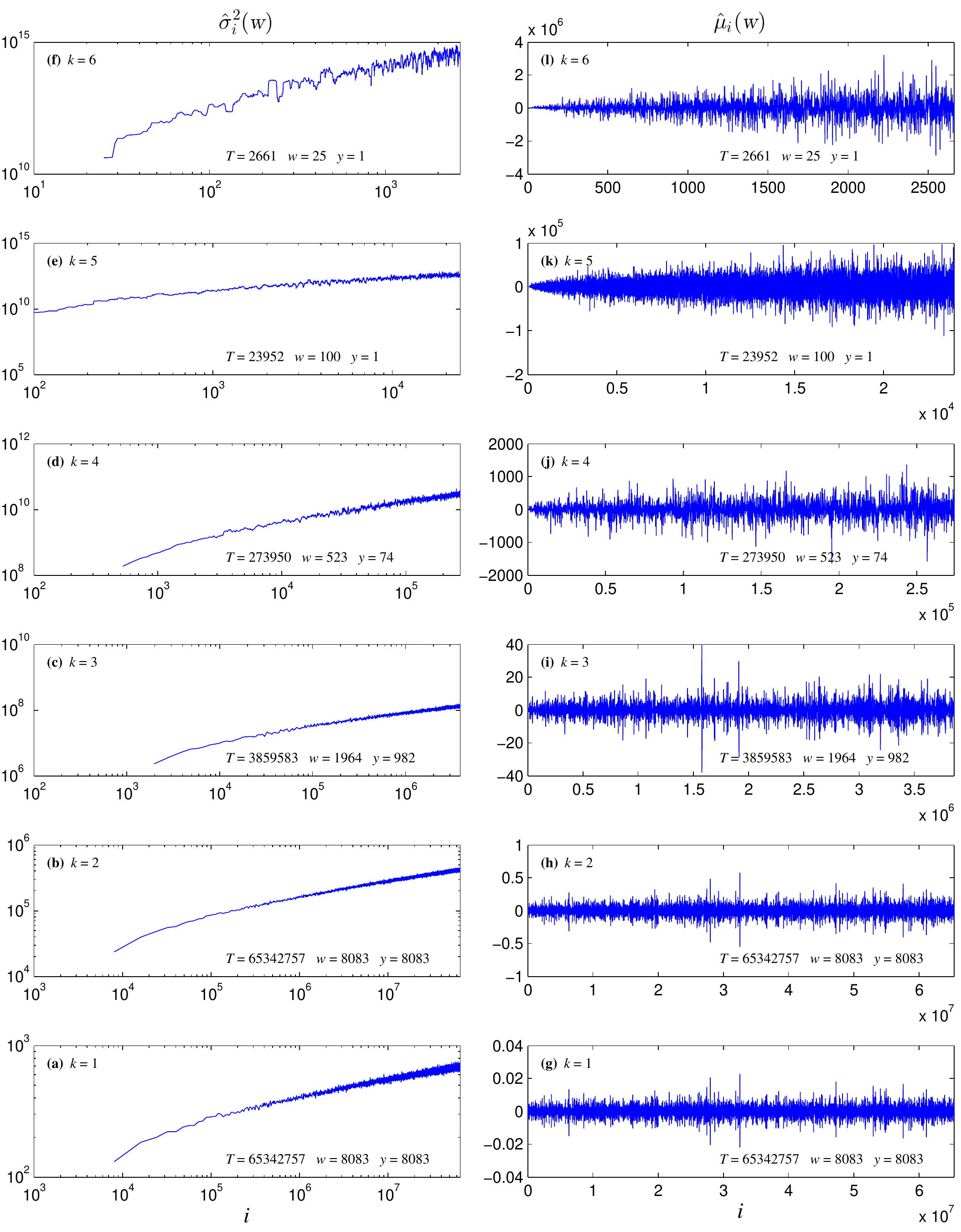}
\caption{Rolling-sample variance and mean of $\daleth_{1,0,i}^{2,k}$ for $k=1..6$, \textit{T} = sample size, \textit{w} = window width, \textit{y} = step size}
\label{fig:VarMeanDQ_K1_6_v02}
\end{figure}

Figure \ref{fig:VarMeanDQ_K1_6_v02} shows 12 plots of $\hat{\sigma}_i^2(w)$ and $\hat{\mu}_i(w)$ for $\daleth_{1,0,i}^{2,k}$. Subplots (a)-(f) show $\hat{\sigma}_i^2(w)$ with $k=1..6$, respectively, and subplots (g)-(l) show $\hat{\mu}_i(w)$ with $k=1..6$, respectively. In each subplot, values for $k$, $T$, $w$, and $y$ are given.  In subplots (a)-(f), $\hat{\sigma}_i(w)$ appears to be non-stationary with a monomial asymptote. In subplots (g)-(l), $\hat{\mu}_i(w)$ appears to be stationary, oscillating around $0$ (which is consistent with the observations of Figure \ref{fig:DELTA_Q_001}). The similarities in fine detail observed in Figure \ref{fig:DELTA_Q_001}(g)-(l) are again seen in Figure \ref{fig:VarMeanDQ_K1_6_v02}(g)-(h), but not in Figure \ref{fig:VarMeanDQ_K1_6_v02}(i)-(l); this is due to varying the values of the parameters $T$, $w$, and $y$ in these subplots.

\subsection[short]{Linear regression.} The phenomena (self-similarity following prime-index order, $k$) exhibited in Table \ref{tab:Daleth_i_k}, Figure \ref{fig:DELTA_Q_001} and Figure \ref{fig:VarMeanDQ_K1_6_v02}(g)-(h) indicates a possible correlation in our data, and we use a simple linear regression as an initial tool for analysis of the relationships between the distributions of the sequences of $\daleth_{1,0,i}^{2,k}$. 

\begin{figure}[hp]
\centering
\includegraphics[width=1\linewidth]{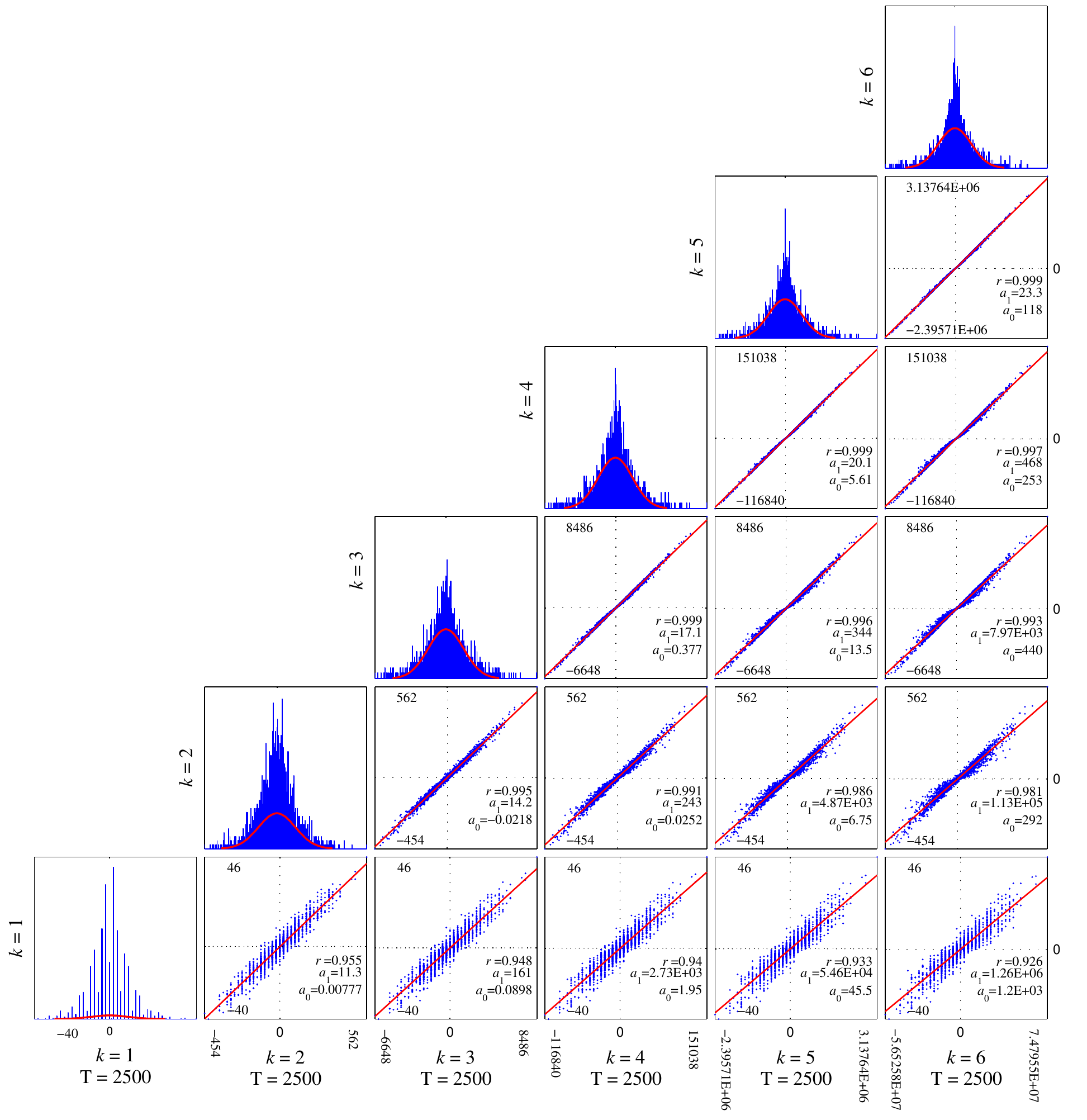}
\caption{Scatterplot matrix and histograms for the pairwise combinations of $\daleth_{1,0,i}^{2,k}$ with $k=1..6$; $T$ = 2500 samples; linear regression lines and Gaussian fits are shown in red; Pearson correlation coefficients $r$, and first-order linear-fit coefficients ($a_0$, $a_1$) are shown in each subplot}
\label{fig:Scatterplots_BIG_dtiff_001}
\end{figure}

Figure \ref{fig:Scatterplots_BIG_dtiff_001} is a scatterplot matrix showing all 15 pairwise combinations of sequences of $\daleth_{1,0,i}^{2,k}$ for $k=1..6$. Histograms, showing the distribution of each sequence $\daleth_{1,0,i}^{2,1},\dots,\daleth_{1,0,i}^{2,6}$, are given along the main diagonal of the matrix. The minima and maxima of the ranges for each pair $(\daleth_{1,0,i}^{2,a},\daleth_{1,0,i}^{2,b})$ are given in the corresponding scatterplot. Each scatterplot also shows the linear regression fit (red line), linear fit coefficients ($a_0$, $a_1$), and Pearson's correlation $r$.

The close linear fits and high correlation coefficients seen in the scatterplots of Figure \ref{fig:Scatterplots_BIG_dtiff_001} indicate that a strong positive correlation may exist in all of the pairwise combination of sequences. The weakest correlation, $r = 0.926$, is seen in the bottom-rightmost scatterplot $(\daleth_{1,0,i}^{2,1},\daleth_{1,0,i}^{2,6})$, while the strongest correlation, $r = 0.999$, is seen in the top-rightmost scatterplot $(\daleth_{1,0,i}^{2,5},\daleth_{1,0,i}^{2,6})$. The correlation strength appears to exhibit the following trends: $r$ monotonically decreases along each row in the direction of increasing $k$; $r$ monotonically decreases along each column in the direction of decreasing $k$; and $r$ monotonically increases along each diagonal in the direction of increasing $k$. Due to the strong correlations, it's reasonable to hypothesize that these trends can be extrapolated to larger samples and higher $k$, continuing ad infinitum.

\subsubsection[short]{Modeling the probability distribution.} In Figure \ref{fig:Scatterplots_BIG_dtiff_001}, each histogram subplot is overlayed with a normal distribution fit (red curve). The normal model is initially chosen since the distributions in all six histograms appear unimodal and symmetric. However, the sharp peaks and fat tails of the distribution data, compared to the normal fits, indicate that the distribution of $\daleth_{1,0,i}^{2,k}$ is leptokurtic over our sample, hence, the normal distribution model is not the ideal choice for a fit. 

In searching for a better fit, Holdom\cite{Holdom-Scale-Invariant-2009}, Dahmen et al.\cite{Dahmen-Similarity-In-The-Statistics-2001}, Kumar\cite{Kumar-Information-Entropy-2003}, Wolf et al.\cite{Wolf1996,Wolf-Nearest_Neighbor-2014}, Scafetta et al.\cite{Scafetta2004} and others\cite{Timberlake-Is-There-Quantum-2008,Selvam-Universal-2014} have observed that the histograms of prime difference functions, e.g., $\Delta_1^1(p_i)$ and $\Delta_1^2(p_i)$, exhibit exponential or Poisson probability distributions. We find that this behavior applies also to the finite differences of PIPs. Hence, an improved model for the distribution of our example case $\daleth_{1,0,i}^{2,k}$ can be obtained with the Laplace (or double exponential) distribution $P(x)=\frac{1}{2b}e^{-|x-\mu|/b}$, where $P(x)$ is the probability density function (PDF). 

\begin{figure}[hp]
\psfrag{n}{T}
\centering
\includegraphics[width=1\linewidth]{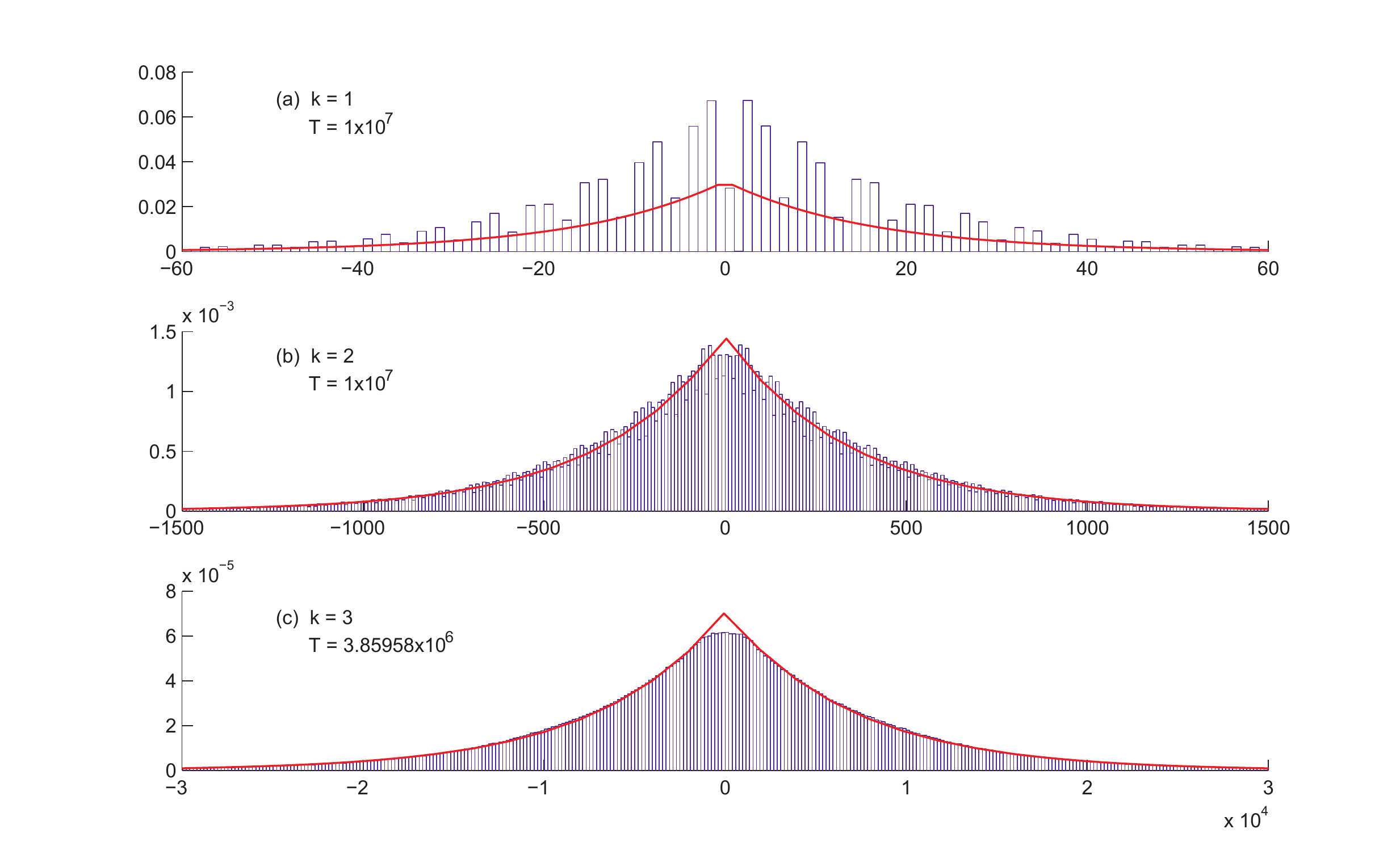}

\includegraphics[width=1\linewidth]{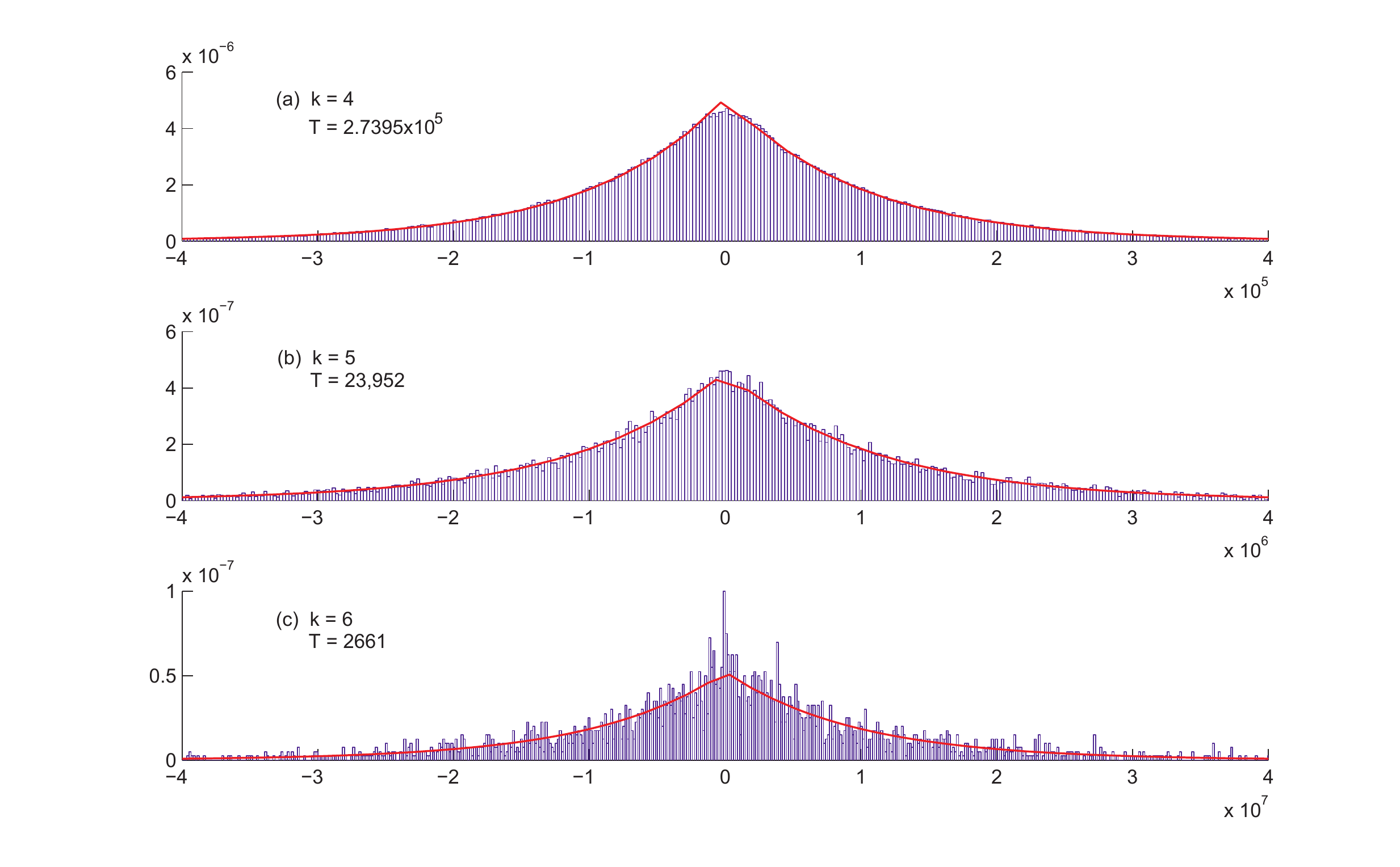}
\caption{PDF histograms (blue bars) and  Laplace distribution curve fits (red lines) for $\daleth_{1,0,i}^{2,k}$ with $k=1..6$ and $T$ = sample size}
\label{fig:LaplaceFits}
\end{figure}

Figure \ref{fig:LaplaceFits} shows PDF histograms of $\daleth_{1,0,i}^{2,k}$ with $k=1..6$. In each histogram, $T$ is the sample size and the domain is $i=1..T$. Due to the large numbers of bins, the outer wings of the histograms have been excluded in order to help visualize the fine detail of central regions of the distributions. In each of the six cases, the Laplace distribution shows a good fit with the histograms of $\daleth_{1,0,i}^{2,k}$. 

\subsubsection[short]{\label{Periodicity in prob dists}Periodicity in the probability distribution.} Periodic (or ``oscillatory'') behavior in the probability distribution histograms of the first and second finite differences of primes, i.e., $\Delta_1^1(p_n)$ and $\Delta_1^2(p_n)$, is an active area of research.\cite{Ares-Hidden-2006,Dahmen-Similarity-In-The-Statistics-2001,Szpiro-TheGaps-2004,Szpiro-Peaks-and_Gaps-2007,Wolf1996,Wolf-Nearest_Neighbor-2014,Kumar-Information-Entropy-2003,Bershadskii-Hidden-2011}\footnote{See footnote \ref{Liang2006}.} Specifically, spikes in the histograms of prime differences are seen to occur when $\Delta_1^1(p_n)\equiv{}0\Mod{3}$ (i.e., oscillation at period 3). Likewise, dips in the histograms of prime increments (differences of differences) are seen to occur when $\Delta_1^2(p_n)\equiv{}0\Mod{6}$ (i.e., oscillation at period 6). Ares and Castro \cite{Ares-Hidden-2006} show that this behavior directly results from the combination of Dirichlet's Theorem \footnote{Given an arithmetic progression of terms $an+b$, for $n=1, 2,\dots{}$, the series contains an infinite number of primes if $a$ and $b$ are relatively prime} and local structure (i.e., that for every prime $p>3$, $p\equiv{}\pm{}1\Mod{6}$).

\begin{figure}[h]
\includegraphics[width=1\linewidth]{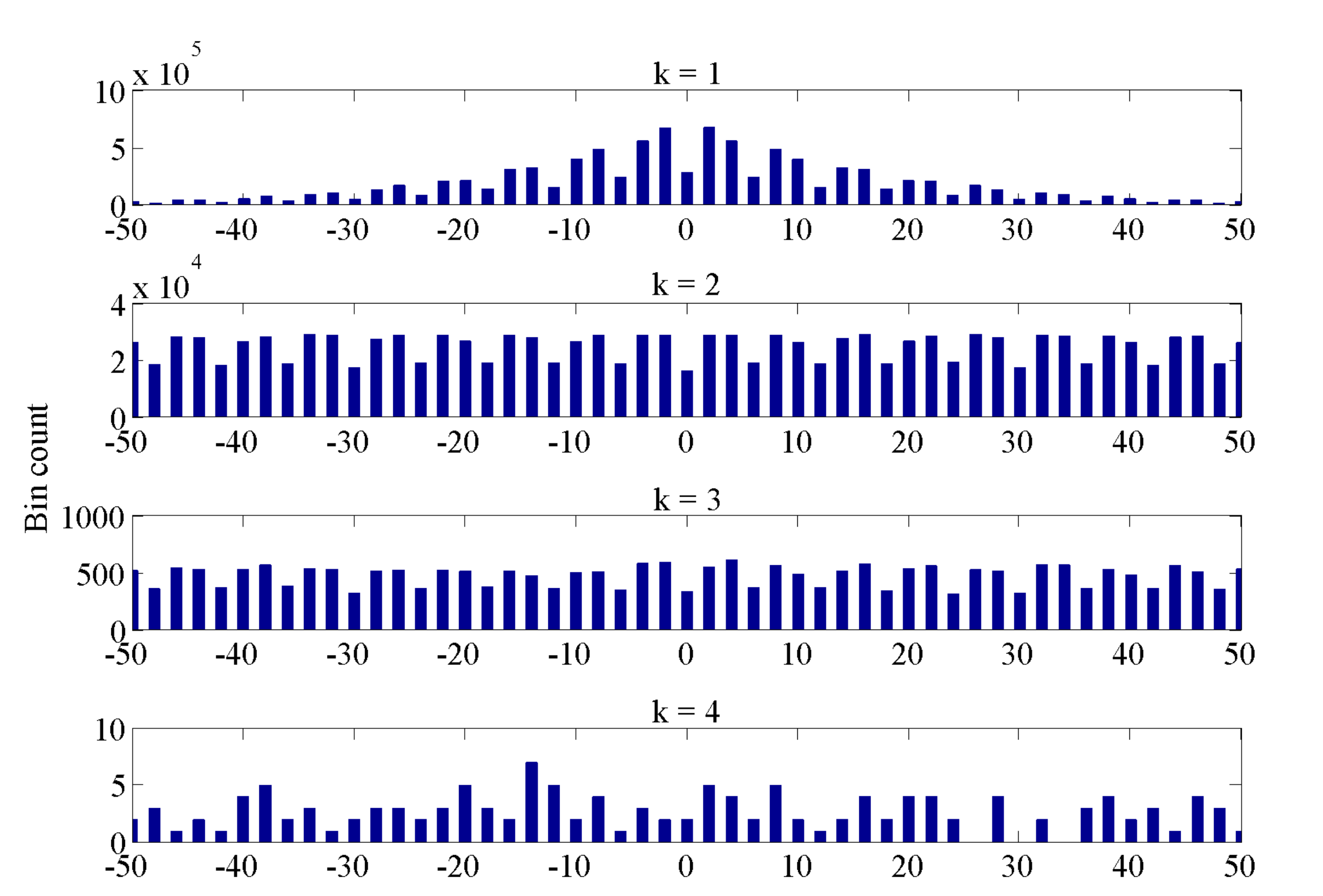}
\caption{Detailed views of the probability distribution histograms for $\daleth_{1,0,i}^{2,k}$ with $k=1..4$; period-6 oscillatory behavior is visible in subplots $k=1-3$, while subplot $k=4$ has insufficient bin counts to display the periodicity in the distribution}
\label{fig:Bars1_4ShowMod}
\end{figure}

Ares and Castro \cite{Ares-Hidden-2006} further conjecture that this periodic behavior, ``should hold also for non-consecutive primes\dots'', and we find that this also appears to be the case for second-differenced PIPs, $\daleth_{1,0,i}^{2,k}$. Figure \ref{fig:Bars1_4ShowMod} shows probability density histograms for $\daleth_{1,0,i}^{2,k}$ with $k=1..4$. Sample sizes $T$ are identical to those given in Figure \ref{fig:LaplaceFits} for the corresponding $k$ values. In order to highlight the fine detail in the histograms, in all four subplots the bin width is set to 1 and bin numbers are limited to a range of $\daleth_{1,0,i}^{2,k}=-50..50$. In subplots $k=1,2,3$, periodic behavior is clearly visible with dips occurring at $0\Mod{6}$. For these three subplots, the average number of counts in bins $-50..50$ is roughly: $5\times{}10^5$ for $k=1$; $2\times{}10^4$ for $k=2$; and $500$ for $k=3$. These counts are sufficient to reveal the fine detail in the probability distribution and its periodic behavior. In contrast, however, for the fourth subplot, $k=4$, the range of $\daleth_{1,0,i}^{2,4}$ is so large that our sample (i.e., the first 1.3 billion primes) provides fewer than $10$ counts in each of bins $-50..50$, thus yielding an insufficient number of counts to display the oscillatory behavior in the distribution. 

\section{Fractal Plots of $\daleth$} 
 
We now investigate the fractal behavior of $\daleth_{hsi}^{nk}$, focusing on the example case $\daleth_{1,0,i}^{2,k}$. We will examine the function both in its raw form as well as after filtering to stationarize the variance.

\subsection[short]{Stationarization of the variance.} It is of interest to test whether the observed correlation among the pairwise combinations of sequences of $\daleth_{1,0,i}^{2,k}$ (see Figure \ref{fig:Scatterplots_BIG_dtiff_001}) is statistically significant or an erroneous result of non-Gaussian distribution and variance nonstationarity (see Figure \ref{fig:VarMeanDQ_K1_6_v02}). Further, in order to obtain useful statistics on $\daleth_{hsi}^{nk}$ (e.g., for forecasting large primes), it will be necessary to detrend the variance $\hat{\sigma}_i^2(w)$, and possibly any higher non-stationary moments, in a given sample. Attempting to model and detrend the moments of $\daleth_{hsi}^{nk}$ is beyond the scope of the present work, however, we can stationarize $\hat{\sigma}_i^2(w)$ and test the correlation by using the sign function,

\[
\hphantom{--------------}\text{sgn}(x)= 
\begin{cases}
    1        & \text{for } x > 0\\
    0        & \text{for } x = 0\\
    -1       & \text{for } x < 0.\\
\end{cases} 
\]

\begin{figure}[hp]
\centering
\includegraphics[width=1\linewidth]{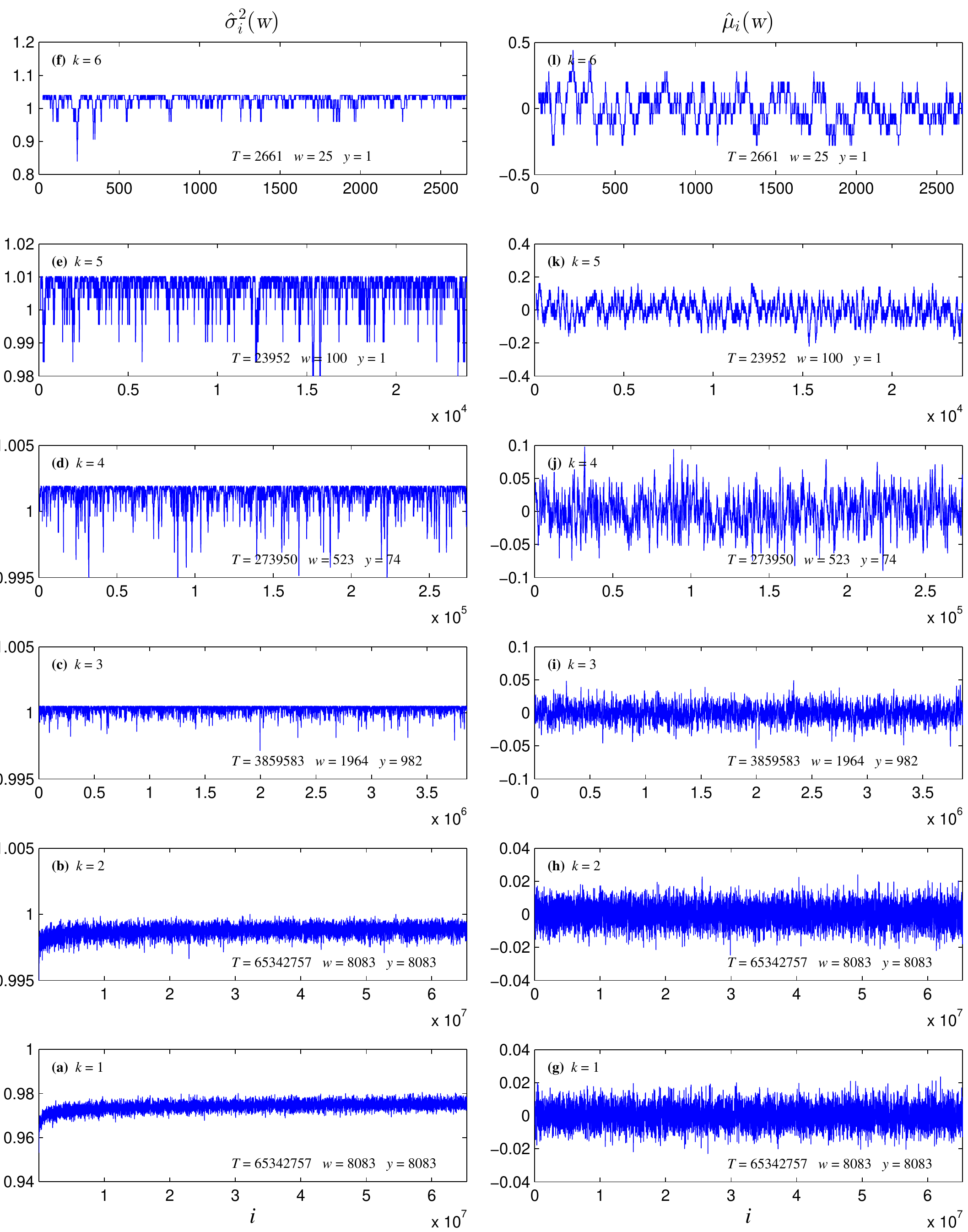}
\caption{Rolling-sample variance and mean of sgn$(\daleth_{1,0,i}^{2,k})$ for $k=1..6$; \textit{T} = sample size, \textit{w} = window width, \textit{y} = step size}
\label{fig:VarMeanDQ_K1_6_SIGNUM_v01}
\end{figure}

Figure \ref{fig:VarMeanDQ_K1_6_SIGNUM_v01} shows 12 plots of moving-average $\hat{\sigma}_i^2(w)$ and $\hat{\mu}_i(w)$ for sgn$(\daleth_{1,0,i}^{2,k})$. All parameters are identical to those given in Figure \ref{fig:VarMeanDQ_K1_6_v02} (which shows $\hat{\sigma}_i^2(w)$ and $\hat{\mu}_i(w)$ for unfiltered $\daleth_{1,0,i}^{2,k}$). It can be seen for $k=1..6$ that the variance is now stationarized at $\hat{\sigma}_i^2(w) \approx{} 1$, while $\hat{\mu}_i(w)$ continues to oscillate around $0$ as previously shown for the unfiltered case.

\subsection[short]{Fractal plot of stationarized $\daleth$.} 
\begin{figure}[h]
\centering
\includegraphics[width=1\linewidth]{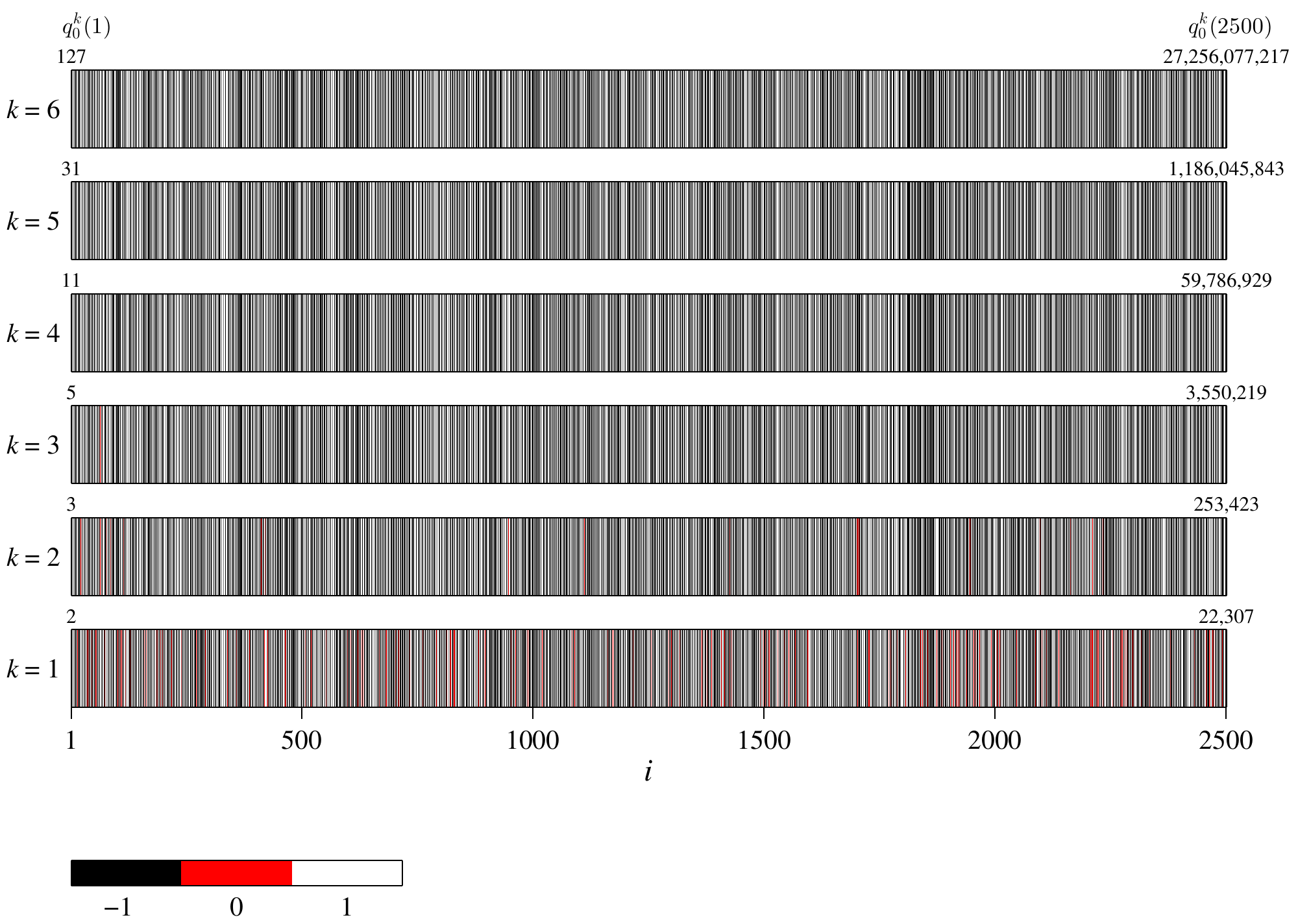}
\caption{Fractal gridplot of sgn$(\daleth_{1,0,i}^{2,k})$ mapped onto a 3-color colormap where $h=1$, $n=2$, $s=0$, $k=1..6$ and $i=1..2500$}
\label{fig:BWFract_002}
\end{figure}

Having stationarized the variance of $\daleth_{1,0,i}^{2,k}$, we now examine its correlation behavior in the context of a fractal plot. Figure \ref{fig:BWFract_002} is a vertical array of six one-dimensional horizontally-oriented grid plots. The data points in the grid plots are stretched vertically into narrow bands in order to give the plots sufficient height for easy viewing. Since the sign function reduced the range of $\daleth_{1,0,i}^{2,k}$ to just three values $(-1,0,1)$, we map sgn$(\daleth_{1,0,i}^{2,k})$ onto a 3-color colormap as follows: sgn$(\daleth_{1,0,i}^{2,k})=1$ (white); sgn$(\daleth_{1,0,i}^{2,k})=0$ (red); and sgn$(\daleth_{1,0,i}^{2,k})=-1$ (black). Each grid plot displays a unique $k$th sequence of $\daleth_{1,0,i}^{2,k}$, and each vertical band (data point) represents a unique $i$th value of $\daleth_{1,0,i}^{2,k}$ of its respective $k$th sequence. In all of the six grid plots, the domain of $\daleth_{1,0,i}^{2,k}$ is $i = 1..2500$ and the $k$ values are labeled to the left of the plots. At the top of each grid plot, on the left and right ends, the values of $q_{0,1}^k=q_0^k(1)$ and $q_{0,2500}^k=q_0^k(2500)$, are marked. Since $\daleth_{1,0,i}^{2,k} = \Delta_1^2(q_{0,i}^k)$, these markings are provided to help convey the scale and range of $q_{0,i}^k$ over our sample. For example, on the bottom grid plot, $k=1$ and $q_{0,i}^1$ ranges from 2 to 22,307. For the top grid plot, $k=6$ and $q_{0,i}^6$ ranges from 127 to 27,256,077,217.

In Figure \ref{fig:BWFract_002}, the quasi-self-similarity of $\daleth_{1,0,i}^{2,k}$, and scale-invariance following prime-index order $k$, are unmistakable. While the upper limit of $q_{0,i}^k$ spans more than six orders of magnitude (i.e., from  $q_{0,2500}^1\approx{}2.2\times{}10^4$ to $q_{0,2500}^6\approx{}2.7\times{}10^{10}$), all six 2500-element sequences are almost identical. In fact, the similarity between adjacent grid plots appears to improve as $k$ gets larger; this is consistent with the monotonically increasing correlation $r$ values observed along the matrix diagonals of Figure \ref{fig:Scatterplots_BIG_dtiff_001}. And, unlike typical fractals, scale invariance of the structure does not follow a power law, but instead follows \mbox{prime-index order $k$.}

\subsection[short]{Fractal plot of unfiltered $\daleth$.} The quasi-self-similarity seen in Figure \ref{fig:BWFract_002} (in which the variance of $\daleth_{1,0,i}^{2,k}$ was stationarized) suggests that the strong correlation observed in Figure \ref{fig:Scatterplots_BIG_dtiff_001} is valid, and not an artifact of a non-Gaussian probability distribution or nonstationarity of the moments of $\daleth_{1,0,i}^{2,k}$. Therefore, it is of interest to apply the format used in Figure \ref{fig:BWFract_002} to the unfiltered (variance non-stationary) $\daleth_{1,0,i}^{2,k}$ data in order to examine the range of the function in finer detail.

\begin{figure}[Ht]
\includegraphics[width=1\linewidth]{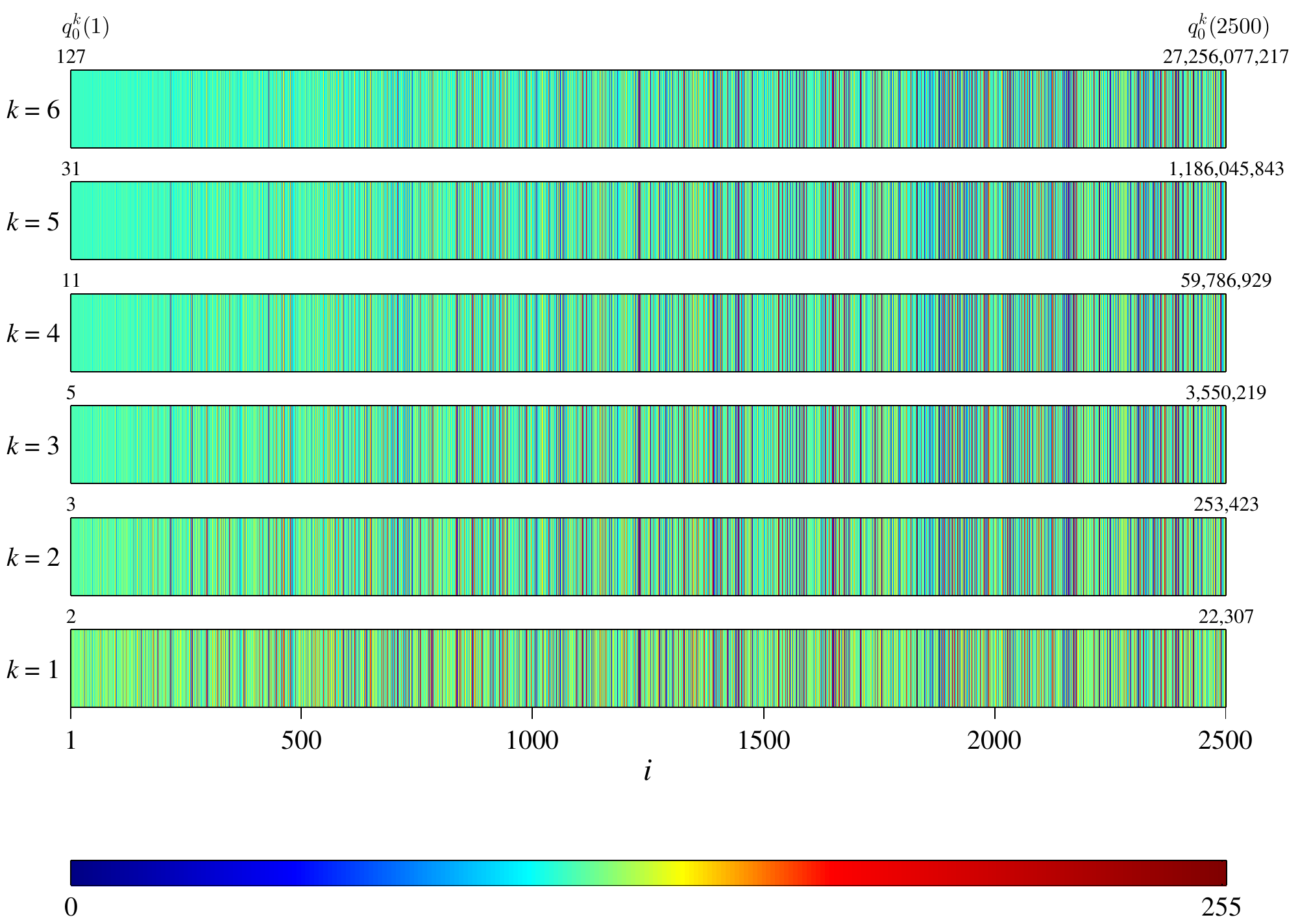}
\caption{Fractal gridplot of $\daleth_{1,0,i}^{2,k}$ mapped onto an 8-bit colormap where $h=1$, $n=2$, $s=0$, $k=1..6$ and $i=1..2500$}
\label{fig:ColorFract_009}
\end{figure}

Over our sample, $\daleth_{1,0,i}^{2,k}$ ranges from 0 to approximately $\pm{}8\times{}10^7$ (see Figure \ref{fig:DELTA_Q_001}(f)), a span of more than eight orders of magnitude. In Figure \ref{fig:BWFract_002}, this range was reduced to just three values (-1,0,1) by the sign function. We now increase the resolution of the range of $\daleth_{1,0,i}^{2,k}$ from three to 256 levels, and map the data to an 8-bit colormap. We use a modified version of the ``Jet'' 8-bit colormap provided in Matlab\textsuperscript{\textregistered} software, which runs from blue to red and includes cyan, green, yellow and orange. In Figure \ref{fig:ColorFract_009}, all parameters are identical to those of Figure \ref{fig:BWFract_002}, except that the range of each of the $k$th $\daleth_{1,0,i}^{2,k}$ sequences is now downsampled and scaled to the interval [0..255] as follows,

\begin{equation}\label{Daleth 256}
\hphantom{---------}\nonumber\daleth_{1,0,i\text{ (256 levels)}}^{2,k}
=\text{nint}\left(\frac{\daleth_{1,0,i}^{2,k}-\text{min}\left(\daleth_{1,0,i}^{2,k}\right)}
{\text{max}\left(\daleth_{1,0,i}^{2,k}\right)-\text{min}\left(\daleth_{1,0,i}^{2,k}\right)}
\cdot{}255\right),
\end{equation}

where $k=1..6$; $i=1..2500$ and nint($x$) is the nearest-integer function.  

\pagebreak
In Figure \ref{fig:ColorFract_009}, the self-similarity and scale invariance of $\daleth_{1,0,i}^{2,k}$ are again indisputable. Although color-bands representing approximately all 256 levels are visible in each grid plot, the color is dominated by the green-yellow portion of the spectrum. The left side of each grid plot is particularly washed-out in comparison to the right side; this is a result of the sharp peak of the Laplace probability distribution of $\daleth_{1,0,i}^{2,k}$ being mapped onto the green portion of the 8-bit colormap and the nonstationarity of the variance of $\daleth_{1,0,i}^{2,k}$.

\section{Discussion}\label{Section:Discussion}

\subsection[short]{Zeros of $\daleth$.} Evaluation of the scaling exponents and fractal dimensionality of the distribution of primes is currently an active area of research;\cite{Cattani-Fractal-2010,Holdom-Scale-Invariant-2009,Scafetta2004,Wolf-Nearest_Neighbor-2014}\footnote{See footnote \ref{Selvam2008}.}$^,$\footnote{See footnote \ref{Cloitre2011}.} a result of the recent studies of statistical structure discussed in \mbox{Sections \ref{Section:Introduction}} and \ref{Periodicity in prob dists}. Extending this fractal analysis to the finite differences of prime-indexed primes $\daleth_{hsi}^{nk}$ will require careful treatment of the function's properties, such as quasi-self-similarity, variance nonstationarity, exponential distribution, and scaling by prime-index order, and is beyond the scope of the present work. Nonetheless, we have illustrated some properties of $\daleth_{hsi}^{nk}$ and $q_{si}^k$ that might be useful in addressing this question. In Section 2, we iterated the findings of Bayless et al. \cite{Bayless2013}, on the asymptotic lower bound for $q_{0,i}^2$, to prime-index order $k=6$. In Section 3, the $r$ correlations of pairwise combinations of $\daleth_{1,0,i}^{2,k}$ were observed to increase for pairs of sequences having consecutive values of $k$. 

Another potentially useful measure for quantifying the fractality of $\daleth_{hsi}^{nk}$ is its distribution of zeros; that is, the values of $i$ and $k$ for which $\daleth_{hsi}^{nk}=0$. For our example case of $h=1$, $s=0$, and $n=2$, the zeros of $\daleth_{1,0,i}^{2,k}$ occur when $2q_{0,i+1}^k=q_{0,i+2}^k+q_{0,i}^k$. This form of recurrence relation has a historical background in the ``balanced primes''\cite{Broughan-Shifted-2012,Ghusayni-subsets-of-primes-2012}. Balanced primes are those primes $p_n$ for which $2p_n = p_{n+j}+p_{n-j}$, and it is conjectured that the number of balanced primes is infinite.\footnote{Wikipedia, Balanced prime, 2014, \url{http://en.wikipedia.org/w/index.php?title=Balanced_prime&oldid=603643701}, [Online; accessed  12-April-2014]}


Erd\"{o}s and Pomerance conjectured the following about balanced primes: ``Is it true that for $n>n_0$ there always is an $i$ for which $2p_n=p_{n+i}+p_{n-i}$? The answer is almost certainly affirmative.''\cite{Erdos-Some-Problems-1971} 

\begin{figure}[h]
\centering
\includegraphics[width=1\linewidth]{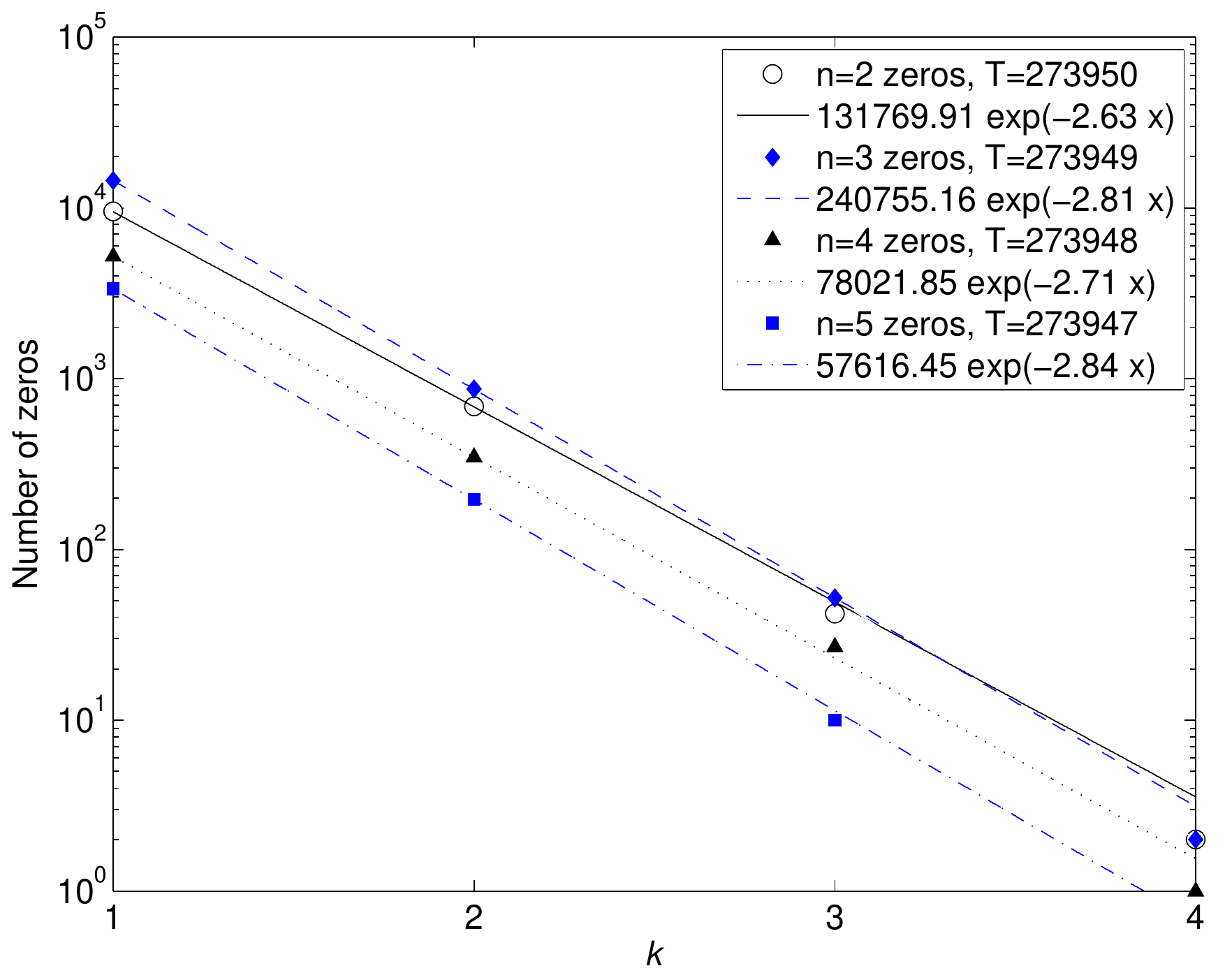}
\caption{Density of zeros of $\daleth_{1,0,i}^{nk}$ vs. $k$ for $n=2..5$, $T$ samples, and exponential fits}
\label{fig:CountZerosExpFit}
\end{figure}

In examining the distribution of zeros of our example case $\daleth_{1,0,i}^{2,k}$, we are interested in their density and infinitude, and ask: for every $k\geq{}1$ is there an infinite number of $i$'s for which $2q_{0,i+1}^k=q_{0,i+2}^k+q_{0,i}^k$? In Figure \ref{fig:BWFract_002}, the zeros of $\daleth_{1,0,i}^{2,k}$ are highlighted with red bands; numerous zeros can be seen in the $k=1$ gridplot with density rapidly diminishing for each sequence of increasing $k$. Figure \ref{fig:CountZerosExpFit} is a plot of the density of zeros of $\daleth_{1,0,i}^{nk}$ vs. $k$ for $k=1..4$. Plots are given for orders of finite differences $n=2..5$. Exponential fits to the data and samples $T$ are shown in the legend. The $R^2$'s of the fits are all $1.0000$. Although our sample only yields four data points for $n=2,3,4$ and three data points for $n=5$, the strong exponential fits give no reason to expect that the number of zeros is finite for any order of $k$ or $n$.

\clearpage

\subsection[short]{Globality.} Up to this point, we have mostly focused on a local case of $\daleth_{hsi}^{nk}$ with parameters fixed at $h=1$, $s=0$, and $n=2$. The prime-index order $k$ has been the chief parameter examined, as it is the fractal scaling parameter of $\daleth_{1,0,i}^{2,k}$, while $i$ is the index of the function.

Returning to the topic of local structure discussed in Section \ref{Section:Introduction}, consider that for our example case $\daleth_{1,0,i}^{2,k}$ and any $k\geq1$, the primes $p_n\in{}q_{0,i}^k$ encode `special' information (i.e., the distribution of $\daleth_{1,0,i}^{2,1}$), while all other primes $p_m \notin{}q_{0,i}^k$ apparently do not. This can easily be tested by replacing the primes of a given $q_{0,i}^k$ sequence with nearby primes $p_m \notin{}q_{0,i}^k$; the result is that the correlations of $(\daleth_{1,0,i}^{2,k},\daleth_{1,0,i}^{2,b\neq{}k})$ will be decimated. Additionally, in Figures \ref{fig:DELTA_Q_001}(g)-\ref{fig:DELTA_Q_001}(l), it can be seen that the values of $\daleth_{1,0,i}^{2,k}$ oscillate rapidly about zero, indicating the sensitivity of $\daleth_{1,0,i}^{2,k}$ to $q_{0,i}^k$. Thus, we may say that $\daleth_{1,0,i}^{2,k}$ is a structure with \textit{local information at the prime-index order (or $k$th) place}.  

Now, we are interested in building a global picture of $\daleth_{hsi}^{nk}$ and will examine its behavior when the parameters $s, n$ and $h$ are freely varied.

\subsubsection[short]{Shift invariance.} As shown in \eqref{eq:Cases}, the underlying index set of any $k$th sequence $q_{si}^{k\geq{2}}$ is $p_i$ (the ordered set of primes). Now, let us perform a left-shift operation $S^{\ast{}}$ on $p_i$ as follows:

\begin{singlespace}
\begin{eqnarray}
\hphantom{--------}\nonumber{}S^{\ast{}0}(p_1,p_2,p_3,\dots)&\rightarrow&(p_1,p_2,p_3,\dots)=(2,3,5,\dots)\\
\nonumber{}S^{\ast{}1}(p_1,p_2,p_3,\dots)&\rightarrow&(p_2,p_3,p_4,\dots)=(3,5,7,\dots)\\
\nonumber{}&&\vdots{}\\
\nonumber{}S^{\ast{}s}(p_1,p_2,p_3,\dots)&\rightarrow&(p_{s+1},p_{s+2},p_{s+3},\dots).
\end{eqnarray}
\end{singlespace}

Then, define a left-shift operation on the PIPs, $q_i^k$. For $k=2$, let

\begin{singlespace}
\begin{align*}
\nonumber{}S^{\ast{}0}(q_i^2)=q_{0,i}^2=(q_{0,1}^2,q_{0,2}^2,q_{0,3}^2,\dots)
   &\rightarrow(p_{p_1},p_{p_2},p_{p_3},\dots)=(3,5,11,\dots)\\
\nonumber{}S^{\ast{}1}(q_i^2)=q_{1,i}^2=(q_{1,1}^2,q_{1,2}^2,q_{1,3}^2,\dots)
   &\rightarrow(p_{1+p_{2}},p_{1+p_{3}},p_{1+p_{4}},\dots)=(7,13,19\dots{})\\
\nonumber{}&\vdots{}\\
\nonumber{}S^{\ast{}s}(q_i^2)=q_{si}^2=(q_{s,1}^2,q_{s,2}^2,q_{s,3}^2,\dots)
   &\rightarrow(p_{s+p_{s+1}},p_{s+p_{s+2}},p_{s+p_{s+3}},\dots).
\end{align*}
\end{singlespace}

In general, define the $k$th-order left-shifted PIPs as

\begin{singlespace}
\begin{equation}\label{eq:q_i cases shifted}
 \hphantom{----------}S^{\ast{}s}(q_i^k)=q_{si}^k\rightarrow
\begin{cases}
    p_{s+i},             & \text{if } k = 1,\\
    p_{s+p_{s+i}},         & \text{if } k = 2,\\
    \vdots          & \vdots\\
    p_{s+p_{s+_{\ddots_{p_{s+i}}}}},   & \text{if } k = m.                 
\end{cases} 
\end{equation}
\end{singlespace}

Table \ref{tab:index-shifted q^k} shows an array of the values of $q_{si}^2$ for $k=2$, $s=0..11$ and a domain of $i=1..30$. The $s=0$, $i=1..20$ elements of Table \ref{tab:index-shifted q^k} are identical to the values in the $k=2$ column of Table \ref{tab:q_k^i}, as both of these data correspond to $(q_{0,1}^2,q_{0,2}^2,\dots,q_{0,20}^2)$. In the other columns, the sequences of primes shown in Tables \ref{tab:q_k^i} and \ref{tab:index-shifted q^k} are unique and suggest that, with appropriate choice of $s$, any prime $p_{s+i}$ can be a member of a corresponding sequence of index-set shifted PIPs.

\begin{table}[ph]
\centering
\fontsize{10pt}{0}
\caption{Values of $q_{si}^2$ for index-set shifts of $s=0..11$}
$\begin{array}{c|cccccccccccc} i & s=0 & s=1 & s=2 & s=3 & s=4 & s=5 & s=6 & s=7 & s=8 & s=9 & s=10 & s=11\\ 
\hline
 1 & 3 & 7 & 17 & 29 & 47 & 61 & 83 & 101 & 127 & 163 & 179 & 223\\ 2 & 5 & 13 & 23 & 43 & 59 & 79 & 97 & 113 & 157 & 173 & 211 & 239\\ 3 & 11 & 19 & 41 & 53 & 73 & 89 & 109 & 151 & 167 & 199 & 233 & 251\\ 4 & 17 & 37 & 47 & 71 & 83 & 107 & 149 & 163 & 197 & 229 & 241 & 271\\ 5 & 31 & 43 & 67 & 79 & 103 & 139 & 157 & 193 & 227 & 239 & 269 & 311\\ 6 & 41 & 61 & 73 & 101 & 137 & 151 & 191 & 223 & 233 & 263 & 307 & 349\\ 7 & 59 & 71 & 97 & 131 & 149 & 181 & 211 & 229 & 257 & 293 & 347 & 359\\ 8 & 67 & 89 & 127 & 139 & 179 & 199 & 227 & 251 & 283 & 337 & 353 & 397\\ 9 & 83 & 113 & 137 & 173 & 197 & 223 & 241 & 281 & 331 & 349 & 389 & 421\\ 10 & 109 & 131 & 167 & 193 & 211 & 239 & 277 & 317 & 347 & 383 & 419 & 433\\ 11 & 127 & 163 & 191 & 199 & 233 & 271 & 313 & 337 & 379 & 409 & 431 & 463\\ 12 & 157 & 181 & 197 & 229 & 269 & 311 & 331 & 373 & 401 & 421 & 461 & 491\\ 13 & 179 & 193 & 227 & 263 & 307 & 317 & 367 & 397 & 419 & 457 & 487 & 541\\ 14 & 191 & 223 & 257 & 293 & 313 & 359 & 389 & 409 & 449 & 479 & 523 & 593\\ 15 & 211 & 251 & 283 & 311 & 353 & 383 & 401 & 443 & 467 & 521 & 587 & 613\\ 16 & 241 & 281 & 307 & 349 & 379 & 397 & 439 & 463 & 509 & 577 & 607 & 619\\ 17 & 277 & 293 & 347 & 373 & 389 & 433 & 461 & 503 & 571 & 601 & 617 & 647\\ 18 & 283 & 337 & 367 & 383 & 431 & 457 & 499 & 569 & 599 & 613 & 643 & 659\\ 19 & 331 & 359 & 379 & 421 & 449 & 491 & 563 & 593 & 607 & 641 & 653 & 683\\ 20 & 353 & 373 & 419 & 443 & 487 & 557 & 587 & 601 & 631 & 647 & 677 & 787\\ 21 & 367 & 409 & 439 & 479 & 547 & 577 & 599 & 619 & 643 & 673 & 773 & 821\\ 22 & 401 & 433 & 467 & 541 & 571 & 593 & 617 & 641 & 661 & 769 & 811 & 857\\ 23 & 431 & 463 & 523 & 569 & 587 & 613 & 631 & 659 & 761 & 809 & 853 & 863\\ 24 & 461 & 521 & 563 & 577 & 607 & 619 & 653 & 757 & 797 & 839 & 859 & 941\\ 25 & 509 & 557 & 571 & 601 & 617 & 647 & 751 & 787 & 829 & 857 & 937 & 953\\ 26 & 547 & 569 & 599 & 613 & 643 & 743 & 773 & 827 & 853 & 929 & 947 & 997\\ 27 & 563 & 593 & 607 & 641 & 739 & 769 & 823 & 839 & 919 & 941 & 991 & 1033\\ 28 & 587 & 601 & 631 & 733 & 761 & 821 & 829 & 911 & 937 & 983 & 1031 & 1061\\ 29 & 599 & 619 & 727 & 757 & 811 & 827 & 907 & 929 & 977 & 1021 & 1051 & 1097\\ 30 & 617 & 719 & 751 & 809 & 823 & 887 & 919 & 971 & 1019 & 1049 & 1093 & 1151\end{array}$
\label{tab:index-shifted q^k}
\end{table}

\begin{figure}
\centering
\includegraphics[width=1\linewidth]{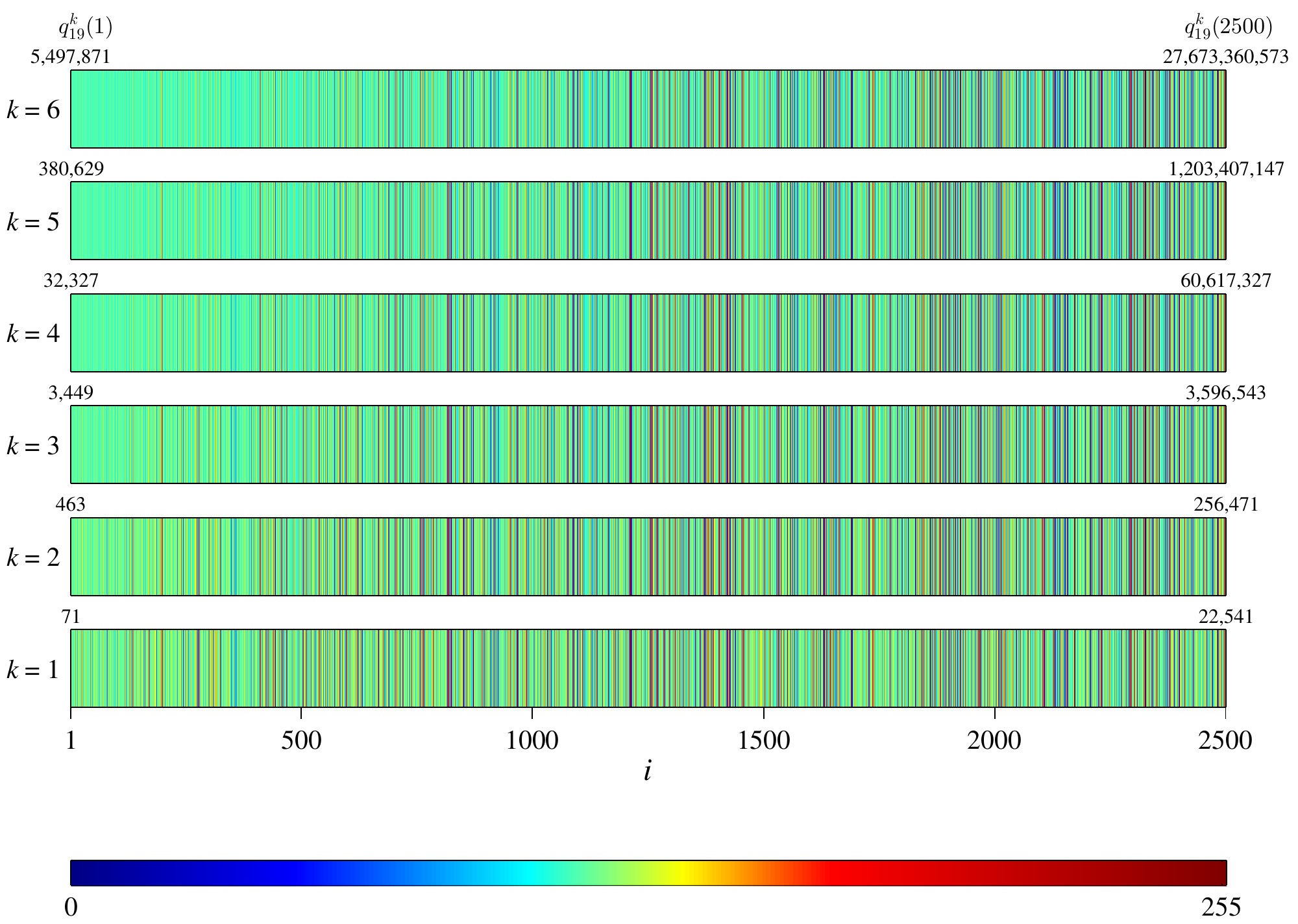}
\caption{Fractal gridplot of $\daleth_{1,19,i}^{2,k}$ mapped onto an 8-bit colormap where $h=1$, $n=2$, $s=19$, $k=1..6$, and $i=1..2500$}
\label{fig:ColorFract_D2h1s19}
\end{figure}

\pagebreak
Using the definition of \eqref{eq:q_i cases shifted}, we now examine $\daleth_{hsi}^{nk}$ for different values of $s$. Figure \ref{fig:ColorFract_D2h1s19} is a fractal plot of $\daleth_{1,19,i}^{2,k}$ with an index-set shift of $s=19$, $h=1$, $n=2$, $k=1..6$, and a domain of $i=1..2500$. The same 8-bit colormap of Figure \ref{fig:ColorFract_009} is used here and in the following figures to emphasize the relative changes in the distributions and moments of $\daleth_{hsi}^{nk}$ as its parameters are varied. It is instructive to compare Figure \ref{fig:ColorFract_D2h1s19} (showing $\daleth_{1,19,i}^{2,k}$ with $s=19$) to Figure \ref{fig:ColorFract_009} (showing $\daleth_{1,0,i}^{2,k}$ with $s=0$). Even though the ranges of the $k$th sequences in these two figures are vastly different, the only significant change is that all six grid plots in Figure \ref{fig:ColorFract_D2h1s19} are slightly shifted to the left by 19. 

\begin{figure}
\centering
\includegraphics[width=1\linewidth]{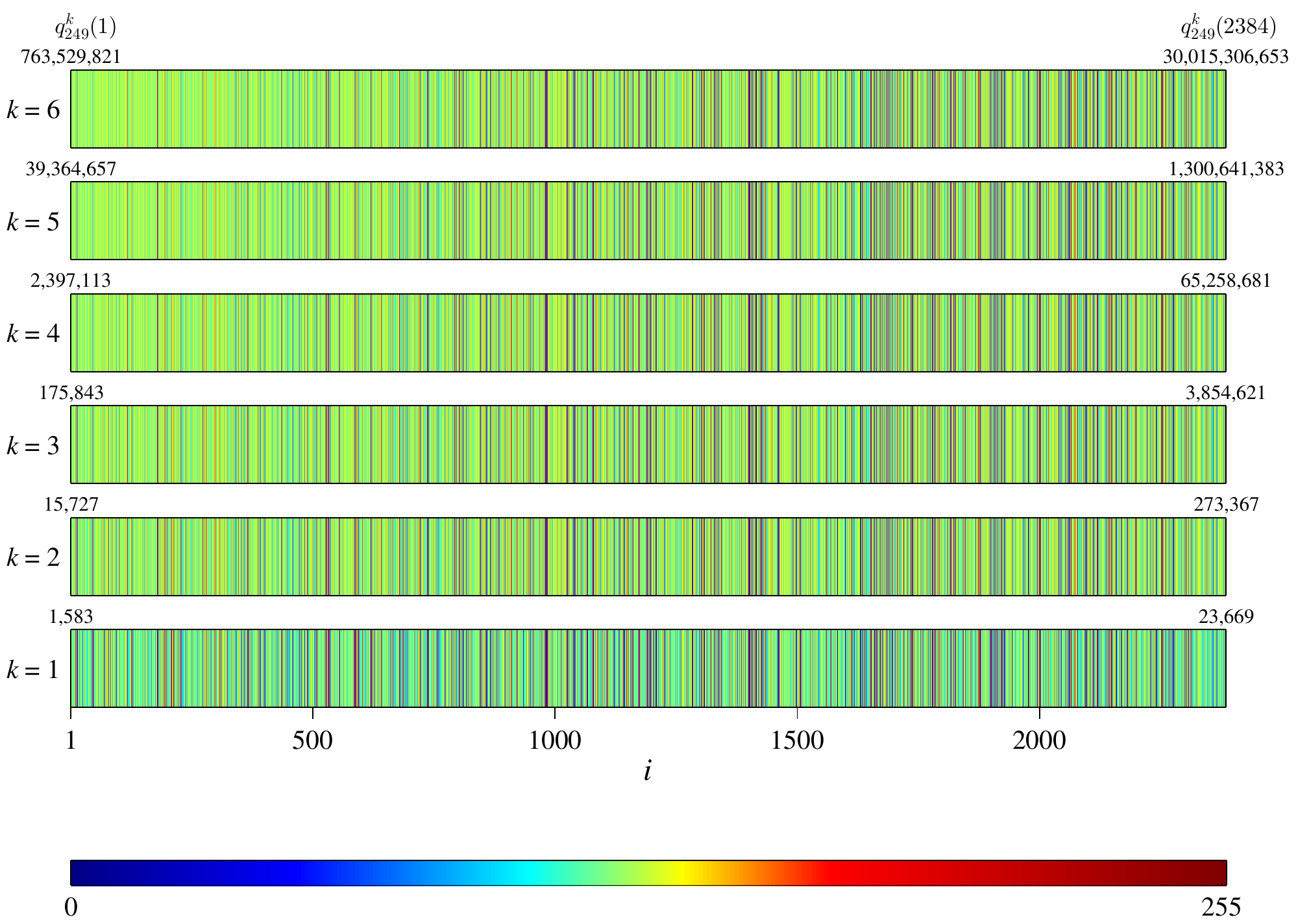}
\caption{Fractal gridplot of $\daleth_{1,249,i}^{2,k}$ mapped onto an 8-bit colormap where $h=1$, $n=2$, $s=249$, $k=1..6$, and $i=1..2384$}
\label{fig:ColorFract_D2h1s249}
\end{figure}

Figure \ref{fig:ColorFract_D2h1s249}, shows $\daleth_{1,249,i}^{2,k}$ where the index-set shift is now increased to $s=249$. The parameters $h$, $n$, and $k$ are identical to those in Figures \ref{fig:ColorFract_009} and \ref{fig:ColorFract_D2h1s19}. In Figure \ref{fig:ColorFract_D2h1s249}, the domain is reduced to $i=1..2384$ due to the limit of our sample at $k=6$. Again, the corresponding grid plots of Figures \ref{fig:ColorFract_D2h1s249} and \ref{fig:ColorFract_009} are similar, with the only major difference being that Figure \ref{fig:ColorFract_D2h1s249} is now shifted to the left by 249. From these observations, we hypothesize that $\daleth_{hsi}^{nk}$ is shift (or translation) invariant on its index set. 

\subsubsection[short]{Spacing and order of finite differences.}
\begin{figure}
\centering
\includegraphics[width=1\linewidth]{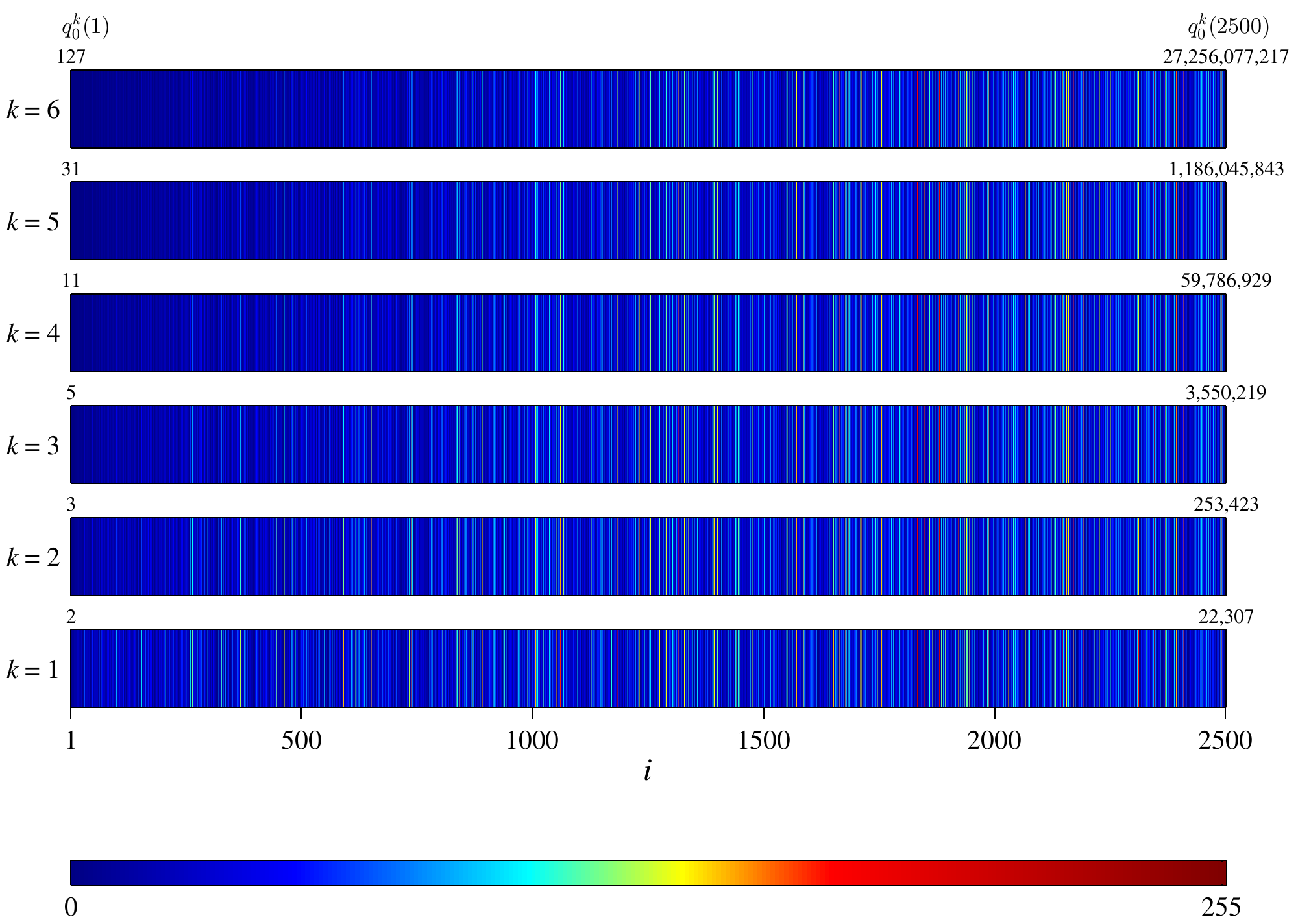}
\caption{Fractal gridplot of $\daleth_{1,0,i}^{1,k}$ mapped onto an 8-bit colormap where $h=1$, $n=1$, $s=0$, $k=1..6$ and $i=1..2500$}
\label{fig:ColorFract_D1h1s0}
\end{figure}

\begin{figure}[p]
\centering
\includegraphics[width=1\linewidth]{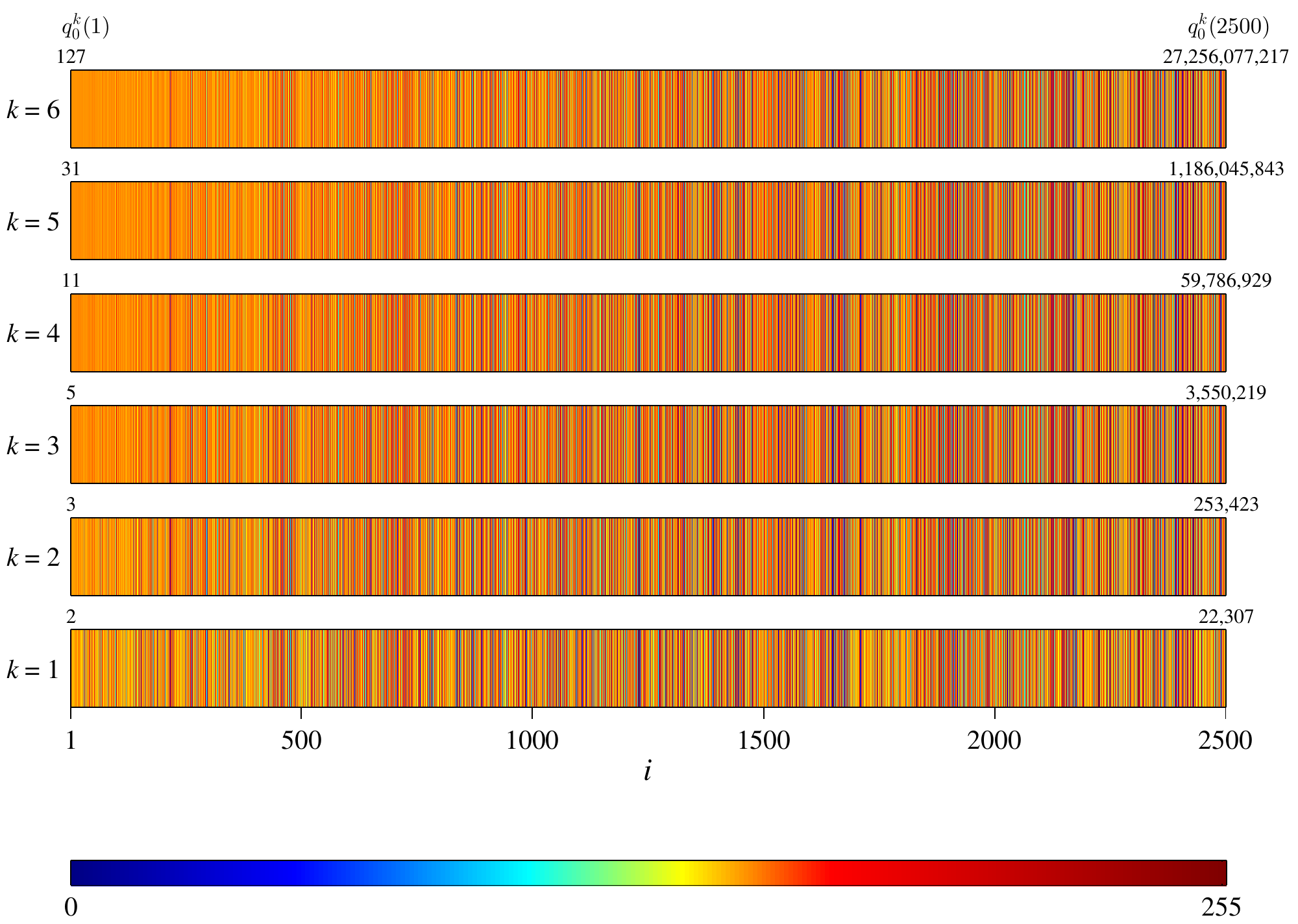}
\caption{Fractal gridplot of $\daleth_{1,0,i}^{3,k}$ mapped onto an 8-bit colormap where $h=1$, $n=3$, $s=0$, $k=1..6$ and $i=1..2500$}
\label{fig:ColorFract_D3h1s0}
\end{figure}

\begin{figure}
\centering
\includegraphics[width=1\linewidth]{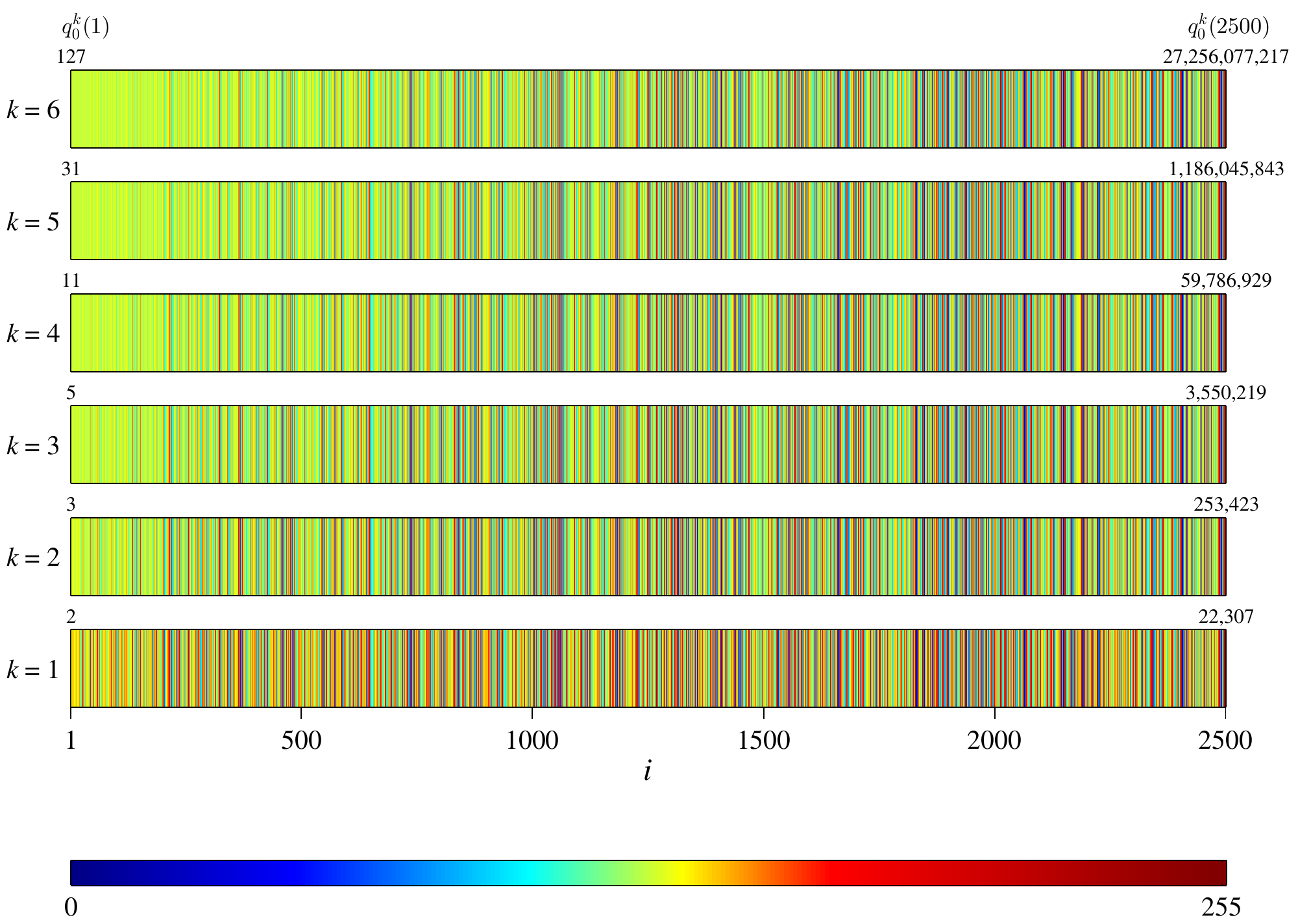}
\caption{Fractal gridplot of $\daleth_{3,0,i}^{2,k}$ mapped onto an 8-bit colormap where $h=3$, $n=2$, $s=0$, $k=1..6$ and $i=1..2500$}
\label{fig:ColorFract_D2h3s0}
\end{figure}

\begin{figure}
\centering
\includegraphics[width=1\linewidth]{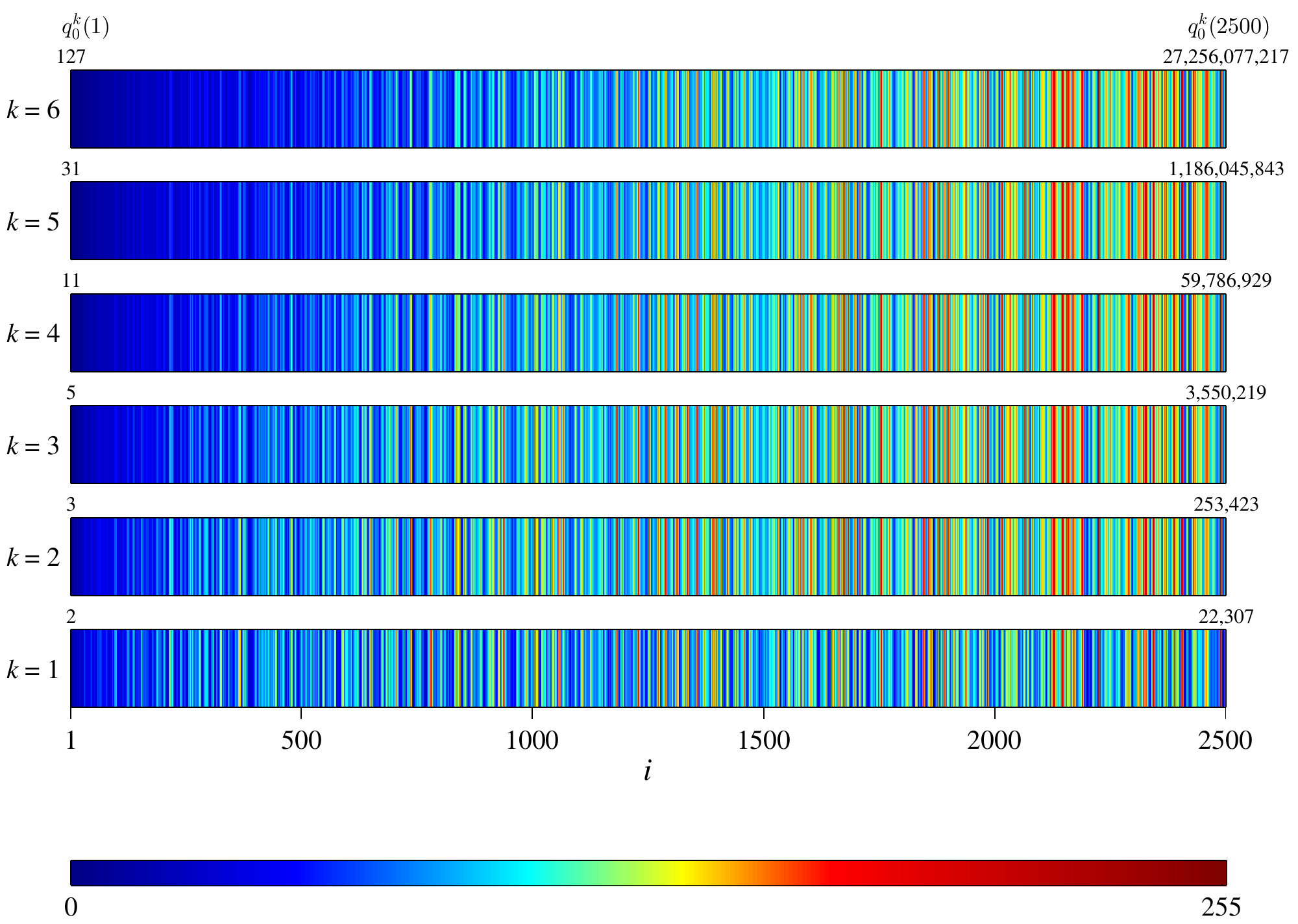}
\caption{Fractal gridplot of $\daleth_{5,0,i}^{1,k}$ mapped onto an 8-bit colormap where $h=5$, $n=1$, $s=0$, $k=1..6$ and $i=1..2500$}
\label{fig:ColorFract_D1h5s0}
\end{figure}

So far, we have kept the spacing $h$ and order $n$ parameters of the finite difference operator (see \eqref{Defn daleth}) fixed at $h=1$ and $n=2$. As shown in Figures \ref{fig:ColorFract_D1h1s0}-\ref{fig:ColorFract_D1h5s0}, we now test the fractal behavior of $\daleth_{hsi}^{nk}$ when these parameters are allowed to vary. In all four figures, $s=0$, and $k=1..6$. The values of $h$ and $n$ in these figures are: $h=1$, $n=1$ (Figure \ref{fig:ColorFract_D1h1s0}); $h=1$, $n=3$ (Figure \ref{fig:ColorFract_D3h1s0}); $h=3$, $n=2$ (Figure \ref{fig:ColorFract_D2h3s0}); and $h=5$, $n=1$ (Figure \ref{fig:ColorFract_D1h5s0}). In Figures \ref{fig:ColorFract_D2h1s19} and \ref{fig:ColorFract_D2h1s249}, it was shown that changes in $s$ result in a translation of the distribution of $\daleth_{hsi}^{nk}$ along the index $i$ axis. Now, in Figures \ref{fig:ColorFract_D1h1s0}-\ref{fig:ColorFract_D1h5s0}, with $s$ fixed and $h$ and $n$ variable, the distribution itself (shown as the range of colors in the grid patterns) widely varies for each configuration of parameters. 

\begin{figure}
\centering
\includegraphics[width=1\linewidth]{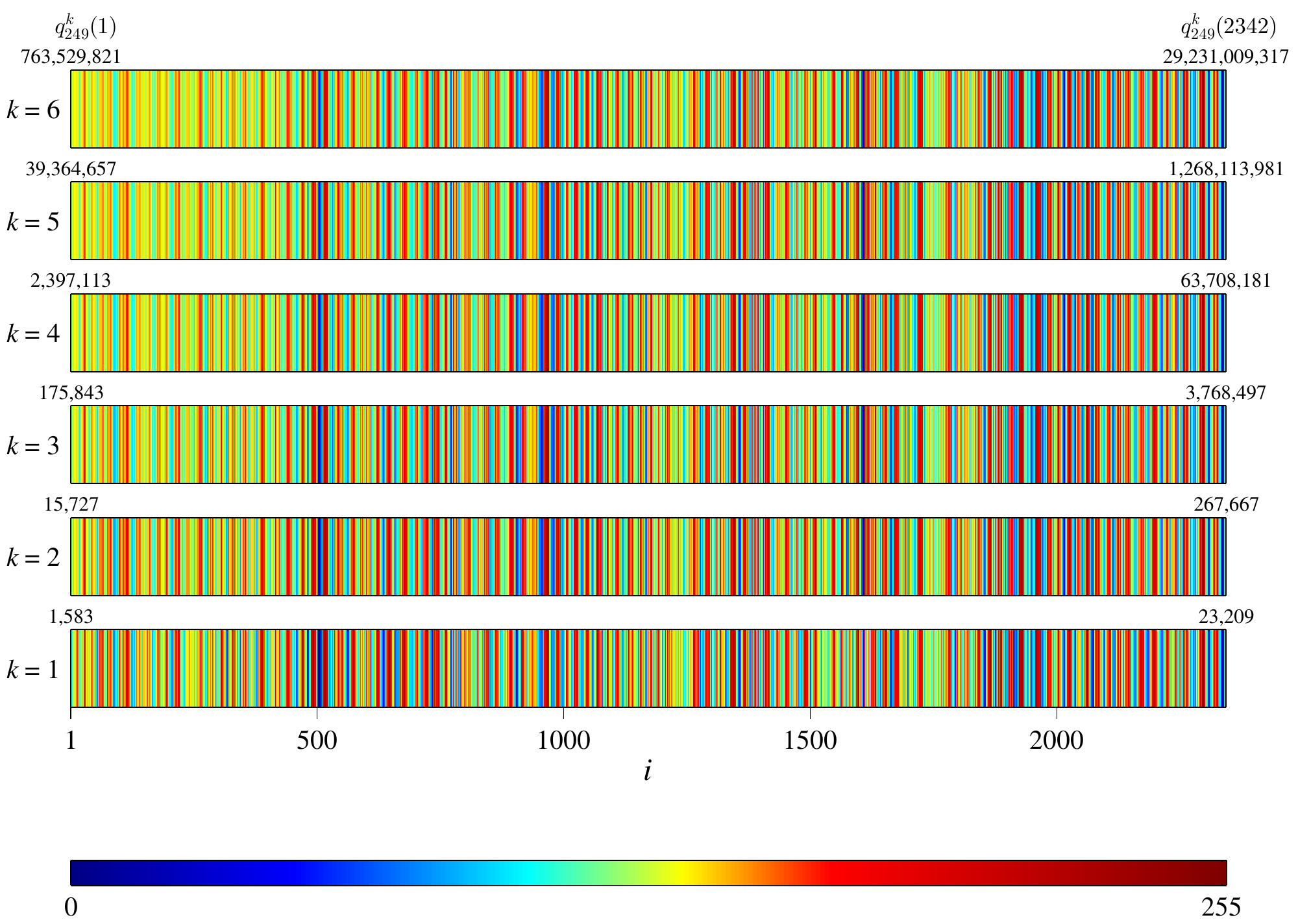}
\caption{Fractal gridplot of $\daleth_{11,249,i}^{4,k}$ mapped onto an 8-bit colormap where $h=11$, $n=4$, $s=249$, $k=1..6$ and $i=1..2342$}
\label{fig:ColorFract_D4h11s249}
\end{figure}

As a final example, Figure \ref{fig:ColorFract_D4h11s249} shows $\daleth_{hsi}^{nk}$ with $h=11$, $n=4$, $s=249$, $k=1..6$, and $i=1..2342$. In spite of now varying all three parameters $h$, $s$ and $n$, the quasi-self-similarity and scaling by prime-index order $k$ are still clearly visible. 

From these examples, we hypothesize that the fractal structure of $\daleth$ extends ad infinitum in $i$ and $k$, and is encoded globally by all primes.

\enlargethispage{10\baselineskip}

\clearpage

\section{Conclusion}
Using a sample of the first 1.3 billion primes, we empirically investigated what appears to be a global quasi-self-similar structure in the distribution of differences of prime-indexed primes, with scaling by prime-index order. We briefly touched on several aspects of the structure to gain a picture of its scope, leaving a plethora more questions yet to be addressed. 

Of the many possible questions, we conclude on the fourth form of structure discussed in Section \ref{Section:Introduction}: if our hypothesis that the fractal behavior of $\daleth$ can be infinitely extended is not proved false, might $\daleth$ then be an exotic structure? Recall that the condition for a structure's being exotic is that it fails to follow the Cram\'{e}r random model. Cram\'er's model predicts that the probability of finding $k$ primes in the interval $[n,n+\text{log}n]$ is $\frac{\lambda{}^k}{k!}\exp{}^{-\lambda}$, which is the Poisson probability distribution\cite{Soundararajan-the-distribution}. However, $\daleth$ is nothing more than a family of sequences, in which each sequence is a linear combination of subsequences of the primes themselves. So, it is no surprise that the sequences of $\daleth$ exhibit exponential or double-exponential distributions, from which the Poisson distribution can be derived.\footnote{Wikipedia, Talk:Poisson distribution, 2013, \url{http://en.wikipedia.org/wiki/Talk:Poisson_distribution, [Online; accessed  13-December-2013]}.} Therefore, $\daleth$ follows the Cram\'{e}r model and does not meet the strict definition of an exotic structure; but ``quasi-exotic structure'' may be a more accurate label, reflecting its pseudorandomness.

The fractality of $\daleth$ may be of relevance in other areas of mathematics and science, such as refining the accuracy of predicting large primes. Likewise, it may be useful in improving the efficiency of factoring large composites, the difficulty of which is critical to RSA public key data encryption. Examples of fractals are found everywhere in nature. So, a fractal structure in the set of primes, which is arguably Nature's most fundamental phenomenon, might help us better understand our world and universe. 

And, while these findings might raise new questions about the primes, perhaps they'll shed light on some old ones too.

\nopagebreak{
\bibliographystyle{amsplain}
\bibliography{PrimesBIBv16}
}

\raggedbottom

\end{document}